\let\myfrac=\frac
\input eplain
\let\frac=\myfrac
\input amstex
\input epsf




\loadeufm \loadmsam \loadmsbm
\message{symbol names}\UseAMSsymbols\message{,}

\font\myfontdefault=cmr10

\font\mytdmchapfont=cmb10 at 14pt
\font\mytdmheadfont=cmb10 at 10pt
\font\mytdmsubheadfont=cmr10

\magnification 1200
\newif\ifinappendices
\newif\ifundefinedreferences
\newif\ifchangedreferences
\newif\ifloadreferences
\newif\ifmakebiblio
\newif\ifmaketdm

\undefinedreferencestrue
\changedreferencesfalse


\loadreferencestrue
\makebibliofalse
\maketdmfalse

\def\headpenalty{-400}     
\def\proclaimpenalty{-200} 

%
%

\def\alphanum#1{\ifcase #1 _\or A\or B\or C\or D\or E\or F\or G\or H\or I\or J\or K\or L\or M\or N\or O\or P\or Q\or R\or S\or T\or U\or V\or W\or X\or Y\or Z\fi}
\def\gobbleeight#1#2#3#4#5#6#7#8{}

\newwrite\references
\newwrite\tdm
\newwrite\biblio

\newcount\chapno
\newcount\headno
\newcount\subheadno
\newcount\procno
\newcount\figno
\newcount\citationno

\def\setcatcodes{%
\catcode`\!=0 \catcode`\\=11}%

\ifloadreferences
    {\catcode`\@=11 \catcode`\_=11%
    \global\def\_@citation@HarveyLawsonI{1}
\global\def\_@citation@CaffNirSprI{2}
\global\def\_@citation@CaffNirSprII{3}
\global\def\_@citation@CaffNirSprV{4}
\global\def\_@citation@GuanSpruckO{5}
\global\def\_@citation@GuanSpruckI{6}
\global\def\_@citation@GuillPoll{7}
\global\def\_@citation@Klingenberg{8}
\global\def\_@citation@LabI{9}
\global\def\_@citation@RosSpruck{10}
\global\def\_@citation@SmiCGC{11}
\global\def\_@citation@SmiNLD{12}
\global\def\_@citation@SmiAAT{13}
\global\def\_@citation@SmiGPP{14}
\global\def\_@citation@Spruck{15}
\global\def\_@citation@TrudWang{16}
\global\def\_@citation@White{17}
\global\def\_@proc@LemmaTrudWangCompactness{1.1}
\global\def\_@proc@ThmExistenceII{1.2}
\global\def\_@head@HeadImmersedSubmanifolds{2}
\global\def\_@proc@CorInvertibility{2.2}
\global\def\_@head@HeadConvexImmersions{3}
\global\def\_@proc@LemmaBoundaryConvexityIsSemiConvexity{3.1}
\global\def\_@proc@LemmaClosednessOfGraphProperty{3.2}
\global\def\_@proc@LemmaOpennessOfGraphProperty{3.3}
\global\def\_@head@HeadConvexityInHigherCodimension{4}
\global\def\_@head@FirstOrderUpperBounds{5}
\global\def\_@proc@LemmaBoundaryCompactnessOfConvexImmersions{5.1}
\global\def\_@proc@PropMPEmbeddedHypersurface{5.2}
\global\def\_@proc@PropAngleBddBelow{5.4}
\global\def\_@proc@PropSecondStep{5.5}
\global\def\_@proc@PropMPNoSelfIntersection{5.6}
\global\def\_@proc@PropFinalStep{5.7}
\global\def\_@proc@PropMPUpperBoundOfInfimum{5.8}
\global\def\_@head@HeadParabolicLimits{6}
\global\def\_@proc@CorTransversalityOfFixedPart{6.1}
\global\def\_@proc@PropMPEnvelopeLiesInsideSigma{6.2}
\global\def\_@proc@PropGeodesicsWithBothEndsInBoundary{6.3}
\global\def\_@head@HeadSemiConvexity{7}
\global\def\_@proc@PropClosednessOfSemiConvexity{7.1}
\global\def\_@proc@LemmaJacobiFieldsLieAbove{7.2}
\global\def\_@head@HeadImmersedBoundaries{8}
\global\def\_@proc@PropDensityOfGenericity{8.1}
\global\def\_@proc@LemmaBoundaryCompactnessWhenNotEmbedded{8.2}
\global\def\_@proc@PropFirstOrderLowerBound{9.1}
\global\def\_@proc@PropConstructionOfBarrier{9.2}
\global\def\_@proc@PropDistanceInGrassmannianA{9.3}
\global\def\_@proc@PropCorrectionTerm{9.4}
\global\def\_@proc@CorGrassmannianDistanceB{9.6}
\global\def\_@proc@PropContinuityOfNormal{9.7}
\global\def\_@fig@FigConfiguration{1}
\global\def\_@proc@LemmaPreCompactness{10.1}
\global\def\_@proc@LemmaRefinedPrecompactness{10.2}
\global\def\_@head@HeadLocalDeformation{11}
\global\def\_@proc@PropTheSetsAreCompactAndConverge{11.1}
\global\def\_@proc@PropSmoothApproximations{11.3}
\global\def\_@head@HeadLocalAndGlobalRigidity{12}
\global\def\_@proc@DefnLocalAndGlobalRigidity{12.1}
\global\def\_@proc@PropPropertiesOfRigidity{12.2}
\global\def\_@proc@LemmaExistence{12.3}
    }%
\else
    \openout\references=references.tex
\fi

\newcount\newchapflag 
\newcount\showpagenumflag 

\global\chapno = -1 
\global\citationno=0
\global\headno = 0
\global\subheadno = 0
\global\procno = 0
\global\figno = 0

\def\resetcounters{%
\global\headno = 0%
\global\subheadno = 0%
\global\procno = 0%
\global\figno = 0%
}

\global\newchapflag=0 
\global\showpagenumflag=0 

\def\chinfo{\ifinappendices\alphanum\chapno\else\the\chapno\fi}%
\def\headinfo{\ifinappendices\alphanum\headno\else\the\headno\fi}%
\def\subheadinfo{\headinfo.\the\subheadno}
\def\procinfo{\headinfo.\the\procno}
\def\figinfo{\the\figno}        
\def\citationinfo{\the\citationno}%
\def\nextheadno{\global\advance\headno by 1 \global\subheadno = 0 \global\procno = 0}
\def\nextsubheadno{\global\advance\subheadno by 1}
\def\nextprocno{\global\advance\procno by 1 \procinfo}
\def\nextfigno{\global\advance\figno by 1 \figinfo}

{\global\let\noe=\noexpand%
%
%
\catcode`\@=11%
\catcode`\_=11%
\setcatcodes%
!global!def!_@@internal@@makeref#1{%
!global!expandafter!def!csname #1ref!endcsname##1{%
!csname _@#1@##1!endcsname%
!expandafter!ifx!csname _@#1@##1!endcsname!relax%
    !write16{#1 ##1 not defined - run saving references}%
    !undefinedreferencestrue%
!fi}}%
!global!def!_@@internal@@makelabel#1{%
!global!expandafter!def!csname #1label!endcsname##1{%
!edef!temptoken{!csname #1info!endcsname}%
!ifloadreferences%
    !expandafter!ifx!csname _@#1@##1!endcsname!relax%
        !write16{#1 ##1 not hitherto defined - rerun saving references}%
        !changedreferencestrue%
    !else%
        !expandafter!ifx!csname _@#1@##1!endcsname!temptoken%
        !else
            !write16{#1 ##1 reference has changed - rerun saving references}%
            !changedreferencestrue%
        !fi%
    !fi%
!else%
    !expandafter!edef!csname _@#1@##1!endcsname{!temptoken}%
    !edef!textoutput{!write!references{\global\def\_@#1@##1{!temptoken}}}%
    !textoutput%
!fi}}%
!global!def!makecounter#1{!_@@internal@@makelabel{#1}!_@@internal@@makeref{#1}}%
!unsetcatcodes%
}
\makecounter{ch}%
\makecounter{head}%
\makecounter{subhead}%
\makecounter{proc}%
\makecounter{fig}%
\makecounter{citation}%
\def\newref#1#2{%
\def\temptext{#2}%
\edef\bibliotextoutput{\expandafter\gobbleeight\meaning\temptext}%
\global\advance\citationno by 1\citationlabel{#1}%
\ifmakebiblio%
    \edef\fileoutput{\write\biblio{\noindent\hbox to 0pt{\hss$[\the\citationno]$}\hskip 0.2em\bibliotextoutput\medskip}}%
    \fileoutput%
\fi}%
\def\cite#1{%
$[\citationref{#1}]$%
\ifmakebiblio%
    \edef\fileoutput{\write\biblio{#1}}%
    \fileoutput%
\fi%
}%
%
%
%

\let\mypar=\par


\def\raggedleft{\leftskip=0pt plus 1fil \parfillskip=0pt}


\font\lettrinefont=cmr10 at 28pt
\def\lettrine #1[#2][#3]#4%
{\hangafter -#1 \hangindent #2
\noindent\hskip -#2 \vtop to 0pt{
\kern #3 \hbox to #2 {\lettrinefont #4\hss}\vss}}

\font\mylettrinefont=cmr10 at 28pt
\def\mylettrine #1[#2][#3][#4]#5%
{\hangafter -#1 \hangindent #2
\noindent\hskip -#2 \vtop to 0pt{
\kern #3 \hbox to #2 {\mylettrinefont #5\hss}\vss}}


\edef\Pagetitle={Blank}

\headline={\hfil\Pagetitle\hfil}

\footline={\hfil\myfontdefault\folio\hfil}

\def\nextoddpage
{
\newpage%
\ifodd\pageno%
\else%
    \global\showpagenumflag = 0%
    \null%
    \vfil%
    \eject%
    \global\showpagenumflag = 1%
\fi%
}


\def\newchap#1#2%
{%
%
%
\global\advance\chapno by 1%
\resetcounters%
%
%
\newpage%
\ifodd\pageno%
\else%
    \global\showpagenumflag = 0%
    \null%
    \vfil%
    \eject%
    \global\showpagenumflag = 1%
\fi%
\global\newchapflag = 1%
\global\showpagenumflag = 1%
%
%
{\font\chapfontA=cmsl10 at 30pt%
\font\chapfontB=cmsl10 at 25pt%
\null\vskip 5cm%
{\chapfontA\raggedleft\hfil%
{%
\ifnum\chapno=0
    \phantom{%
    \ifinappendices%
        Annexe \alphanum\chapno%
    \else%
        \the\chapno%
    \fi}%
\else%
    \ifinappendices%
        Annexe \alphanum\chapno%
    \else%
        \the\chapno%
    \fi%
\fi%
}%
\par}%
\vskip 2cm%
{\chapfontB\raggedleft%
\lineskiplimit=0pt%
\lineskip=0.8ex%
\hfil #1\par}%
\vskip 2cm%
}%
\edef\Pagetitle{#2}%
%
%
\ifmaketdm%
    \def\temp{#2}%
    \def\tempbis{\nobreak}%
    \edef\chaptitle{\expandafter\gobbleeight\meaning\temp}%
    \edef\mynobreak{\expandafter\gobbleeight\meaning\tempbis}%
    \edef\textoutput{\write\tdm{\bigskip{\noexpand\mytdmchapfont\noindent\chinfo\ - \chaptitle\hfill\noexpand\folio}\par\mynobreak}}%
\fi%
\textoutput%
}


\def\newhead#1%
{%
\ifhmode%
    \mypar%
\fi%
\ifnum\headno=0%
\ifinappendices
    \nobreak\vskip -\lastskip%
    \nobreak\vskip .5cm%
\fi
\else%
    \nobreak\vskip -\lastskip%
    \nobreak\vskip .5cm%
\fi%
\nextheadno%
\ifmaketdm%
    \def\temp{#1}%
    \edef\sectiontitle{\expandafter\gobbleeight\meaning\temp}%
    \edef\textoutput{\write\tdm{\noindent{\noexpand\mytdmheadfont\quad\headinfo\ - \sectiontitle\hfill\noexpand\folio}\par}}%
    \textoutput%
\fi%
\font\headfontA=cmbx10 at 14pt%
{\headfontA\noindent\headinfo\ - #1.\hfil}%
\nobreak\vskip .5cm%
}%


\def\newsubhead#1%
{%
\ifhmode%
    \mypar%
\fi%
\ifnum\subheadno=0%
\else%
    \penalty\headpenalty\vskip .4cm%
\fi%
\nextsubheadno%
\ifmaketdm%
    \def\temp{#1}%
    \edef\subsectiontitle{\expandafter\gobbleeight\meaning\temp}%
    \edef\textoutput{\write\tdm{\noindent{\noexpand\mytdmsubheadfont\quad\quad\subheadinfo\ - \subsectiontitle\hfill\noexpand\folio}\par}}%
    \textoutput%
\fi%
\font\subheadfontA=cmsl10 at 12pt
{\subheadfontA\noindent\subheadinfo\ #1.\hfil}%
\nobreak\vskip .25cm %
}%

%
%


\font\mathromanten=cmr10
\font\mathromanseven=cmr7
\font\mathromanfive=cmr5
\newfam\mathromanfam
\textfont\mathromanfam=\mathromanten
\scriptfont\mathromanfam=\mathromanseven
\scriptscriptfont\mathromanfam=\mathromanfive
\def\mathroman{\fam\mathromanfam}


\font\sf=cmss12

\font\sansseriften=cmss10
\font\sansserifseven=cmss7
\font\sansseriffive=cmss5
\newfam\sansseriffam
\textfont\sansseriffam=\sansseriften
\scriptfont\sansseriffam=\sansserifseven
\scriptscriptfont\sansseriffam=\sansseriffive
\def\mathsf{\fam\sansseriffam}


\font\bftwelve=cmb12

\font\boldten=cmb10
\font\boldseven=cmb7
\font\boldfive=cmb5
\newfam\mathboldfam
\textfont\mathboldfam=\boldten
\scriptfont\mathboldfam=\boldseven
\scriptscriptfont\mathboldfam=\boldfive
\def\mathbf{\fam\mathboldfam}


\font\mycmmiten=cmmi10
\font\mycmmiseven=cmmi7
\font\mycmmifive=cmmi5
\newfam\mycmmifam
\textfont\mycmmifam=\mycmmiten
\scriptfont\mycmmifam=\mycmmiseven
\scriptscriptfont\mycmmifam=\mycmmifive

\def\hexa#1{\ifcase #1 0\or 1\or 2\or 3\or 4\or 5\or 6\or 7\or 8\or 9\or A\or B\or C\or D\or E\or F\fi}
\mathchardef\mathi="7\hexa\mycmmifam7B
\mathchardef\mathj="7\hexa\mycmmifam7C


\font\mymsbmten=msbm10 at 8pt
\font\mymsbmseven=msbm7 at 5.6pt
\font\mymsbmfive=msbm5 at 4pt
\newfam\mymsbmfam
\textfont\mymsbmfam=\mymsbmten
\scriptfont\mymsbmfam=\mymsbmseven
\scriptscriptfont\mymsbmfam=\mymsbmfive

\mathchardef\mybeth="7\hexa\mymsbmfam69
\mathchardef\mygimmel="7\hexa\mymsbmfam6A
\mathchardef\mydaleth="7\hexa\mymsbmfam6B


\def\placelabel[#1][#2]#3{{%
\setbox10=\hbox{\raise #2cm \hbox{\hskip #1cm #3}}%
\ht10=0pt%
\dp10=0pt%
\wd10=0pt%
\box10}}%

\def\placefigure#1#2#3{%
\medskip%
\midinsert%
\vbox{\line{\hfil#2\epsfbox{#3}#1\hfil}%
\vskip 0.3cm%
\line{\noindent\hfil\sl Figure \nextfigno\hfil}}%
\medskip%
\endinsert%
}


\newif\ifinproclaim%
\global\inproclaimfalse%
\def\proclaim#1{%
\medskip%
%
%
\bgroup%
\inproclaimtrue%
\setbox10=\vbox\bgroup\leftskip=0.8em\noindent{\bftwelve #1}\sf%
}

\def\endproclaim{%
\egroup%
\setbox11=\vtop{\noindent\vrule height \ht10 depth \dp10 width 0.1em}%
\wd11=0pt%
\setbox12=\hbox{\copy11\kern 0.3em\copy11\kern 0.3em}%
\wd12=0pt%
\setbox13=\hbox{\noindent\box12\box10}%
\noindent\unhbox13%
\egroup%
\medskip\ignorespaces%
}

\def\proclaim#1{%
\medskip%
\bgroup%
\inproclaimtrue%
\noindent{\bftwelve #1}%
\nobreak\medskip%
\sf%
}

\def\endproclaim{%
\mypar\egroup\penalty\proclaimpenalty\medskip\ignorespaces%
}

\def\noskipproclaim#1{%
\medskip%
\bgroup%
\inproclaimtrue%
\noindent{\bf #1}\nobreak\sl%
}

\def\endnoskipproclaim{%
\mypar\egroup\penalty\proclaimpenalty\medskip\ignorespaces%
}


\def\ninn{{n\in\Bbb{N}}}
\def\minn{{m\in\Bbb{N}}}

\def\proof{{\noindent\bf Proof:\ }}

\def\remark{{\noindent\sl Remark:\ }}
\def\example{{\noindent\sl Example:\ }}

\def\mlimsup{\mathop{{\mathroman LimSup}}}

\def\minf{\mathop{{\mathroman Inf}}}
\def\msf#1{{\mathsf #1}}

\def\qed{~$\square$}
\def\munion{\mathop{\cup}}
\def\minter{\mathop{\cap}}
\def\myitem#1{%
    \noindent\hbox to .5cm{\hfill#1\hss}
}

\catcode`\@=11
\def\Eqalign#1{\null\,\vcenter{\openup\jot\m@th\ialign{%
\strut\hfil$\displaystyle{##}$&$\displaystyle{{}##}$\hfil%
&&\quad\strut\hfil$\displaystyle{##}$&$\displaystyle{{}##}$%
\hfil\crcr #1\crcr}}\,}
\catcode`\@=12

\def\makeop#1{%
\global\expandafter\def\csname op#1\endcsname{{\mathroman #1}}}%

\def\makeopsmall#1{%
\global\expandafter\def\csname op#1\endcsname{{\mathroman{\lowercase{#1}}}}}%

\makeopsmall{ArcTan}%
\makeopsmall{ArcCos}%
\makeop{Arg}%
\makeop{Det}%
\makeop{Log}%
\makeop{Re}%
\makeop{Im}%
\makeop{Dim}%
\makeopsmall{Tan}%
\makeop{Ker}%
\makeopsmall{Cos}%
\makeopsmall{Sin}%
\makeop{Exp}%
\makeopsmall{Tanh}%
\makeop{Tr}%
\makeop{End}%
\makeop{Long}%
\makeop{Ch}%
\makeop{Exp}%
\makeop{Eval}%
\makeop{Lift}%
\makeop{Int}%
\makeop{Ext}%
\makeop{Aire}%
\makeop{Im}%
\makeop{Conf}%
\makeop{Exp}%
\makeop{Mod}%
\makeop{Log}%
\makeop{Sgn}%
\makeop{Adj}%
\makeop{Ext}%
\makeop{Supp}%
\makeop{Int}%
\makeop{Ln}%
\makeop{Dist}%
\makeop{Aut}%
\makeop{Id}%
\makeop{GL}%
\makeop{SO}%
\makeop{Homeo}%
\makeop{Vol}%
\makeop{ext}%
\makeop{Ric}%
\makeop{Hess}%
\makeop{Euc}%
\makeop{Isom}%
\makeop{Max}%
\makeop{SW}%
\makeop{SL}%
\makeop{Long}%
\makeop{Fix}%
\makeop{Wind}%
\makeop{Diag}%
\makeop{dVol}%
\makeop{Symm}%
\makeop{Ad}%
\makeop{Diam}%
\makeop{loc}%
\makeopsmall{Sinh}%
\makeopsmall{Cosh}%
\makeop{Len}%
\makeop{Length}%
\makeop{Conv}%
\makeop{Min}%
\makeop{Area}%

\font\mycirclefont=cmsy7
\def\textcircle{{\raise 0.3ex \hbox{\mycirclefont\char'015}}}

\let\emph=\bf

\hyphenation{quasi-con-formal}

%
%

\ifmakebiblio%
    \openout\biblio=biblio.tex %
    {%
        \edef\fileoutput{\write\biblio{\bgroup\leftskip=2em}}%
        \fileoutput
    }%
\fi%

\newref{HarveyLawsonI}{Harvey F. R., Lawson H. B. Jr., Dirichlet Duality and the Nonlinear Dirichlet Problem, {\sl Comm. Pure Appl. Math.} {\bf 62} (2009), no. 3, 396--443}
\newref{CaffNirSprI}{Caffarelli L., Nirenberg L., Spruck J., The Dirichlet problem for nonlinear second-Order elliptic equations. I. Monge Amp\`ere equation, {\sl Comm. Pure Appl. Math.} {\bf 37} (1984), no. 3, 369--402}
\newref{CaffNirSprII}{Caffarelli L., Kohn J. J., Nirenberg L., Spruck J., The Dirichlet problem for nonlinear second-order elliptic equations. II. Complex Monge Amp\`ere, and uniformly elliptic, equations, {\sl Comm. Pure Appl. Math.} {\bf 38} (1985), no. 2, 209--252}
\newref{CaffNirSprV}{Caffarelli L., Nirenberg L., Spruck J., Nonlinear second-order elliptic equations. V. The Dirichlet problem for Weingarten hypersurfaces, {\sl Comm. Pure Appl. Math.} {\bf 41} (1988), no. 1, 47--70}
\newref{GuanSpruckO}{Guan B., Spruck J., Boundary value problems on $\Bbb{S}^n$ for surfaces of constant Gauss curvature, {\sl Ann. of Math.} {\bf 138} (1993), 601--624}
\newref{GuanSpruckI}{Guan B., Spruck J., The existence of hypersurfaces of constant Gauss curvature with prescribed boundary, {\sl J. Differential Geom.} {\bf 62} (2002), no. 2, 259--287}
\newref{GuillPoll}{Guillemin V., Pollack A., {\sl Differential Topology\/}, Prentice-Hall, Englewood Cliffs, N.J., (1974)}
\newref{Klingenberg}{Klingenberg W., {\sl Lectures on closed geodesics\/}, Grundlehren der Mathematischen Wissenschaften, {\bf 230}, Springer-Verlag, Berlin-New York, (1978)}
\newref{LabI}{Labourie F., Un lemme de Morse pour les surfaces convexes (French), {\sl Invent. Math.} {\bf 141} (2000), no. 2, 239--297}
\newref{RosSpruck}{Rosenberg H., Spruck J., On the existence of locally convex hypersurfaces of constant Gauss curvature in hyperbolic space, {\sl J. Differential Geom.} {\bf 40} (1994), no. 2, 379--409}
\newref{SmiCGC}{Smith G., Constant Gaussian Curvature Hypersurfaces in Hadamard Manifold}
\newref{SmiNLD}{Smith G., The Non-Linear Dirichlet Problem in Hadamard Manifolds, arXiv:0908.3590}
\newref{SmiAAT}{Smith G., An Arzela-Ascoli Theorem for immersed submanifolds, {\sl Ann. Fac. Sci. Toulouse Math.}, {\bf 16}, no. 4, (2007), 817--866}
\newref{SmiGPP}{Smith G., The non-linear Plateau problem in non-positively curved manifolds, arXiv:1004.0374}
\newref{Spruck}{Spruck J., Fully nonlinear elliptic equations and applications to geometry, {\sl Proceedings of the International Congress of Mathematicians}, (Z\"urich, 1994), 1145--1152, Birkh\"a user, Basel, 1995.}
\newref{TrudWang}{Trudinger N. S., Wang X., On locally locally convex hypersurfaces with boundary, {\sl J. Reine Angew. Math.} {\bf 551} (2002), 11--32}
\newref{White}{White B., The space of minimal submanifolds for varying Riemannian metrics, {\sl Indiana Univ. Math. J.} {\bf 40}  (1991),  no. 1, 161--200}

\ifmakebiblio%
    {\edef\fileoutput{\write\biblio{\egroup}}%
    \fileoutput}%
\fi%

%
%
%
\document

\myfontdefault
\global\chapno=1
\global\showpagenumflag=1
\def\Pagetitle{}
\null
\vfill
\def\centre{\rightskip=0pt plus 1fil \leftskip=0pt plus 1fil \spaceskip=.3333em \xspaceskip=.5em \parfillskip=0em \parindent=0em}%
\def\textmonth#1{\ifcase#1\or January\or Febuary\or March\or April\or May\or June\or July\or August\or September\or October\or November\or December\fi}
\font\abstracttitlefont=cmr10 at 14pt
{\abstracttitlefont\centre Compactness for immersions of prescribed\break Gaussian curvature II - geometric aspects\par}
\bigskip
{\centre Graham Smith\par}
\bigskip
{\centre 16 February 2010\par}
\bigskip
{\centre Departament de Matem\`atiques,\par
Facultat de Ci\`encies, Edifici C,\par
Universitat Aut\`onoma de Barcelona,\par
08193 Bellaterra,\par
Barcelona,\par
SPAIN\par}
\bigskip
\noindent{\emph Abstract:\ }We develop a compactness result near the boundary for families of locally convex immersions. We also develop a mod 2 degree theory for immersion of constant (and prescribed) Gaussian curvature with prescribed boundary. These are then used to solve the Plateau problem for immersions of constant (and prescribed) Gaussian curvature in general Hadamard manifolds.
%
%
\par
\vfill
\nextoddpage
\global\pageno=1
\myfontdefault
\def\Pagetitle{\sl Compactness for immersions of prescribed Gaussian curvature II}
\newhead{Introduction}
\noindent In this paper, we prove a compactness result for families of locally strictly convex immersions in general Riemannian manifolds. Our motivation is the compactness result for locally convex immersions in $\Bbb{R}^{n+1}$ developed by Guan and Spruck in \cite{GuanSpruckI} and Trudinger and Wang in \cite{TrudWang}, and which may be expressed as follows: let $(\Gamma_n)_\ninn$ be a sequence of compact, codimension $2$ submanifolds of $\Bbb{R}^{n+1}$ which converges to a compact, codimension $2$ submanifold, $\Gamma_0$, say. Let $(\Sigma_n)_\ninn$ be a sequence of compact, locally convex hypersurfaces in $\Bbb{R}^{n+1}$ such that, for all $n$:
$$
\partial\Sigma_n = \Gamma_n.
$$
\noindent For all $n$, let $\opVol(\Sigma_n)$ and $\opDiam(\Sigma_n)$ be the volume and diameter of $\Sigma_n$ respectively. The reasoning of \cite{GuanSpruckI} and \cite{TrudWang} readily yields:
\proclaim{Theorem \nextprocno}
\noindent Suppose that there exists $B>0$ such that, for all $n$:
$$
\opVol(\Sigma_n),\opDiam(\Sigma_n)\leqslant B.
$$
\noindent Then, there exists a $C^{0,1}$, locally convex immersion, $\Sigma_0$ such that:
\medskip
\myitem{(i)} $\partial\Sigma_0=\Gamma_0$; and
\medskip
\myitem{(ii)} $(\Sigma_n)_\ninn$ subconverges towards $\Sigma_0$.
\endproclaim
\proclabel{LemmaTrudWangCompactness}
\noindent The main interest of this result lies in its applications to the study of the Plateau problem. Indeed, after the progress made in the 80's and early 90's (c.f. \cite{CaffNirSprV}, \cite{GuanSpruckO} and \cite{RosSpruck}), Spruck conjectured in \cite{Spruck} that for any $K>0$, for any $k\in]0,K[$, and for any smooth, compact, codimension $2$ curve, $\Gamma\subseteq\Bbb{R}^{n+1}$, bounding a compact, locally strictly convex immersed hypersurface of Gaussian curvature everywhere bounded below by $K$, the Plateau problem can be solved for $\Gamma$ in the sense that there exists a smooth, locally strictly convex, immersed hypersurface, $\Sigma\subseteq\Bbb{R}^{n+1}$, whose boundary is equal to $\Gamma$ and whose Gaussian curvature is constant and equal to $k$. This conjecture was confirmed simultaneously by Guan and Spruck in \cite{GuanSpruckI} and Trudinger and Wang in \cite{TrudWang}. Both papers use the Perron Method in a similar manner, and in each paper some equivalent of Theorem \procref{LemmaTrudWangCompactness} plays a crucial role.
\medskip
\noindent The most natural extension of the results of \cite{GuanSpruckI} and \cite{TrudWang} is then to construct solutions to the Plateau problem in the case where $\Bbb{R}^{n+1}$ is replaced by any Hadamard manifold. We quickly encounter two significant obstacles to the application of the Perron Method to this more general setting. The first is that the Perron Method as described in \cite{GuanSpruckI} and \cite{TrudWang} relies heavily on the existence of a totally geodesic hypersurface passing through every point and normal to any vector at that point, but this is a strong geometric property only possessed by space forms. The second obstacle, similar to the first but different, is that there exists to date no complete extension of Theorem \procref{LemmaTrudWangCompactness} to general manifolds, and that its current proof is only valid in space forms, since it also relies on the existence of the same totally geodesic hypersurfaces described above.
\medskip
\noindent We therefore require completely different techniques, and we proceed as follows: first, in place of the Perron Method, we develop a parametric version of the continuity method arising essentially from the marriage of the approaches of \cite{GuanSpruckI} and \cite{LabI}. This yields a mod $2$ degree theory and, used in conjunction with the appropriate generalisation of Theorem \procref{LemmaTrudWangCompactness}, discussed presently, and the higher order a-priori estimates already developed by the author in \cite{SmiCGC}, this allows us to solve the Plateau problem for Gaussian curvature in general manifolds.
\medskip
\noindent We now state the result, and for ease of presentation, we restrict attention to the case where the ambient manifold is a Hadamard manifold. Thus, let $M^{n+1}$ be an $(n+1)$-dimensional Hadamard manifold. We recall that an {\bf immersed hypersurface} in $M$ is a pair $\Sigma^n=(i,(S^n,\partial S^n))$ where $(S^n,\partial S^n)$ is a smooth, compact, $n$-dimensional manifold with boundary and $i:S\rightarrow M$ is a smooth immersion and $\Sigma$ is said to be {\bf locally convex} if and only if its shape operator is everywhere positive definite. We say that $\Sigma$ has {\bf generic boundary} if and only if for any $p\neq q\in\partial S$ such that $i(p)=i(q)$:
$$
Di\cdot T_p\partial S \neq Di\cdot T_q\partial S.
$$
\noindent In other words, the two tangent spaces of $\partial\Sigma$ at these points do not coincide. Trivially, every smooth, locally convex immersion can be approximated arbitrarily closely by a smooth, locally convex immersion with generic boundary. We prove:
\proclaim{Theorem \nextprocno}
\noindent Let $\hat{\Sigma}^n=(\hat{\mathi},(\hat{S},\partial\hat{S}))$ be a locally strictly convex, immersed hypersurface in $M$ with generic boundary. Let $\phi\in C^\infty(M)$ be a smooth, positive function such that, for every $p\in\hat{\Sigma}$, the Gaussian curvature of $\hat{\Sigma}$ at $p$ is strictly greater than $\phi(p)$. Suppose that there exists a convex set, $K$, with smooth boundary and an open subset $\Omega\subset\partial K$ such that:
\medskip
\myitem{(i)} $\partial\Omega$ is smooth;
\medskip
\myitem{(ii)} $\Omega^c$ has finitely many connected components; and
\medskip
\myitem{(iii)} $\Sigma^n$ is isotopic by locally strictly convex, immersed hypersurfaces to a finite covering of $\Omega$,
\medskip
\noindent then there exists a locally strictly convex, immersed hypersurface $\Sigma^n$ in $M$ such that:
\medskip
\myitem{(a)} $\partial\Sigma = \partial\hat{\Sigma}$;
\medskip
\myitem{(b)} $\Sigma$ is bounded by $\hat{\Sigma}$; and
\medskip
\myitem{(c)} for every point $p\in\Sigma$, the Gaussian curvature of $\Sigma$ at $p$ is equal to $\phi(p)$.
\endproclaim
\proclabel{ThmExistenceII}
\noindent Remarks:
\medskip
\myitem{(a)} the concept of being bounded (condition $(b)$), is described explicitely in Section \headref{HeadConvexImmersions}. Heuristically, $\Sigma$ is bounded by $\hat{\Sigma}$ if and only if $\hat{\Sigma}$ is, more or less, a graph over $\Sigma$;
\medskip
\myitem{(b)} when the ambient manifold is of dimension greater than $3$, immersed hypersurfaces of constant Gaussian curvature typically do not behave well under passage to the limit. We thus do not expect that an approximation argument may be used to relax the condition of genericity along the boundary;
\medskip
\myitem{(c)} if $(\hat{S},\partial\hat{S})$ is diffeomorphic to the closed unit ball in $\Bbb{R}^n$, then Conditions $(i)$, $(ii)$ and $(iii)$ of Theorem \procref{ThmExistenceII} are automatically satisfied;
\medskip
\myitem{(d)} when $n=2$, and thus when the dimension of the ambient manifold is equal to $3$, Theorem \procref{ThmExistenceII} yields a stronger version of Proposition $5.0.3$ of \cite{LabI}, which itself constitutes the analytic core of that paper; and
\medskip
\myitem{(e)} in general manifolds (of arbitrary curvature), the situation is complicated by the possible existence of conjugate points along geodesics. However, most stages of the argument remain more or less intact, and the result may thus be adapted, albeit with stronger hypothesis, to the more general case.
\medskip
\noindent The continuity method as presented here relies no less heavily on a version of Theorem \procref{LemmaTrudWangCompactness} than does the Perron Method it was designed to replace. However, although the complete generalisation remains unproven, we discover that we only really require an analogous result to hold in some neighbourhood of the boundary: this suffices to obtain a-priori estimates which may then be extended across the whole hypersurface using the maximum principal. It is such a boundary compactness result that we prove and which constitutes the hardest and most technical part of this paper, but which is also of independant interest.
\medskip
\noindent The result itself is rather technical, and in order to state it, we first require some notation. Let $M^{n+1}$ be an $(n+1)$-dimensional Riemannian manifold. Let $\Gamma^{n-1}\subseteq M$ be a compact, codimension $2$ immersed submanifold in $M$ which intersects itself generically, as described above. If $p$ is a point in $\Gamma$ and if $\msf{N}_p$ is a unit vector normal to $\Gamma$ at $p$, then we say that $\msf{N}_p$ is a strictly convex normal if and only if the shape operator of $\Gamma$ with respect to $\msf{N}_p$ is strictly positive definite. We then say that $\Gamma$ is strictly convex if and only if, for every $p\in\Gamma$ there exists a strictly convex normal. For example, when $M$ is $3$ dimensional, and $\Gamma$ is thus a closed curve, $\Gamma$ is convex if and only if its geodesic curvature never vanishes. In this case, the set of strictly convex normals over any point constitutes an open subinterval of the circle of unit normal vectors to $\Gamma$ at that point. If, moreover, $\Gamma$ is oriented, which is the case when $\Gamma$ bounds a locally strictly convex hypersurface, we may define a maximal strictly convex normal, $\msf{N}^+$, which is a continuous vector field over $\Gamma$. Let $\Sigma$ be a locally strictly convex immersed hypersurface in $M$ such that $\partial\Sigma=\Gamma$. Suppose, moreover that the orientation of $\Sigma$ is compatible with that of $\Gamma$, and denote the outward pointing unit normal over $\Sigma$ by $\msf{N}_\Sigma$. We obtain:
\proclaim{Lemma \procref{LemmaBoundaryCompactnessOfConvexImmersions}}
\noindent Choose $\theta>0$. There exists $r>0$, which only depends on $M$, $\Gamma$ and $\theta$ such that if the angle between $\msf{N}_\Sigma$ and $\msf{N}^+$ is always greater than $\theta$, then, for all $p\in\Gamma$, there exists a convex subset $K\subseteq B_r(p)$ such that the connected component of $\Sigma\minter B_r(p)$ containing $p$ is embedded and is a subset of $\partial K$.
\endproclaim
\remark Details of notation and conventions are given in Section \headref{HeadConvexityInHigherCodimension}. Note that it is precisely at this stage that the genericity condition on the boundary is required.
\medskip
\noindent This paper is structured as follows:
\medskip
\myitem{(i)} in Sections $2$ to $4$ we introduce the concepts and notation used in the sequel: in Section $2$, we introduce immersed hypersurfaces and describe the Banach manifold of immersed hypersurfaces modulo reparametrisation; in Section $3$ we introduce locally convex immersions and describe the notion of locally convex immersions bounding other locally convex immersions; and in Section $4$ we develope a higher codimension concept of convexity which is required to understand the boundary conditions used in the sequel;
\medskip
\myitem{(ii)} in Sections $5$ to $9$, which together constitute the most innovative part of the paper, we determine first order a-priori bounds near the boundary for generic, locally convex immersions of prescribed curvature: in Section $5$, restricting to the case where the boundary is embedded, and using the notion of ``semi-convexity'', we obtain a compactness result for convex immersions which yields these a-priori bounds but requires various intuitive but technical propositions whose proofs are deferred to the subsequent two sections; in Section $6$, we obtain technical results using the parabolic limit; in Section $7$, we show that the limit of a sequence of semi-convex sets is also semi-convex; in Section $8$, we show how, under a simple modification, the reasoning of Section $5$ may be adapted to the case where the boundary is immersed and generic; and in Section $9$, we obtain first order lower bounds along the boundary which are important in the sequel for the final (technical) step in obtaining second order bounds over the boundary;
\medskip
\myitem{(iii)} in Sections $10$ to $11$, we recall the results of \cite{SmiCGC} to prove a conditional existence result: in Section $10$, we prove compactness of families of immersions of prescribed curvature; and in Section $11$, we show how Sard's Lemma may be used along with compactness to obtain (generically) solutions which interpolate between isotopic data; and
\medskip
\myitem{(iv)} in Section $12$, we prove the existence of isotopies between the given data and other data for which solutions are known to exist, and, using the concepts of local and global rigidity, we prove Theorem \procref{ThmExistenceII}.
\medskip
\noindent This paper was written while the author was working at the Mathematics Department of the Universitat Aut\`onoma de Barcelona, Bellaterra, Spain.
\goodbreak
\newhead{Immersed Submanifolds and Moduli Spaces}
\noindent Let $M^{n+1}$ be a smooth Riemannian manifold. A (smooth, compact) {\bf immersed submanifold} is a pair $\Sigma:=(i,(S,\partial S))$ where:
\headlabel{HeadImmersedSubmanifolds}
\medskip
\myitem{(i)} $(S,\partial S)$ is an oriented, compact, Riemannian manifold with boundary; and
\medskip
\myitem{(ii)} $i:\Sigma\rightarrow M$ is a smooth immersion (i.e. $Di$ is everywhere injective).
\medskip
\remark in the sequel, all submanifolds of $M$ will be (relatively) compact. Likewise, unless stated otherwise, all submanifolds of $M$ will be smooth.
\medskip
\noindent Let $\Sigma=(i,(S,\partial S))$ and $\Sigma' = (i',(S',\partial S'))$ be two immersed hypersurfaces in $M$. We say that $\Sigma$ and $\Sigma'$ are {\bf equivalent} if and only if there exists a diffeomorphism $\phi:(S,\partial S)\rightarrow(S',\partial S')$ such that:
$$
i'\circ\phi = i.
$$
\noindent Let $\opExp$ be the exponential map of $M$. Let $\msf{N}_\Sigma$ be the outward pointing normal vector field over $\Sigma$. We say that $\Sigma'$ is a {\bf graph} over $\Sigma$ if and only if there exists $f\in C_0^\infty(S)$ and a diffeomorphism $\phi:(S,\partial S)\rightarrow (S',\partial S')$ such that:
$$
i'\circ\phi = \opExp(f\msf{N}_\Sigma).
$$
\noindent In particular, $\Sigma$ and $\Sigma'$ are {\bf equivalent} if and only if $\Sigma'$ is a trivial graph over $\Sigma$.
\medskip
\noindent Let $(\Sigma_n)_\ninn = (i_n,(S_n,\partial S_n)),\Sigma_0=(i_0,(S_0,\partial S_0))$ be immersed submanifolds in $M$. We say that $(\Sigma_n)_\ninn$ {\bf converges} to $\Sigma_0$ if and only if there exists $N\geqslant 0$ and, for all $n\geqslant N$ a diffeomorphism $\phi_n:(S_0,\partial S_0)\rightarrow (S_n,\partial S_n)$ such that $(i_n\circ\phi_n)_{n\geqslant N}$ converges to $i_0$ in the $C^\infty$ sense.
\medskip
\noindent Trivially, if $(\Sigma_n)_\ninn$ converges to $\Sigma_0$, then there exists $N\geqslant 0$, and for all $n\geqslant N$ a vector field $X_n\in\Gamma(i_0^*TM)$ and a diffeomorphism $\phi_n:(S_0,\partial S_0)\rightarrow(S_n,\partial S_n)$ such that:
$$
i_n\circ\phi_n = \opExp(X_n).
$$
\noindent Moreover, $(X_n)_{n\geqslant N}$ tends to $0$ in the $C^\infty$ sense. If $\Sigma_n$ and $\Sigma_0$ have the same boundary for all $n$, then, increasing $N$ if necessary, $X_n$ may always be chosen to be normal to $\Sigma_0$ and vanishing along $\partial S_0$. In other words, $\Sigma_n$ is a graph over $\Sigma_0$ for sufficiently large $n$.
\medskip
\noindent Let $(\Gamma_t)_{t\in[0,1]}=(j_t,G_t)_{t\in[0,1]}$ be a smooth family of (exact) immersed submanifolds without boundary in $M$. We denote by $\hat{\Cal{M}}$ the family of all pairs $(t,\Sigma)$ where $t\in I$ and $\Sigma$ is an immersed submanifold in $M$ such that $\partial\Sigma = \Gamma_t$. For all $t\in[0,1]$, let $\hat{\Cal{M}}_t$ be the fibre of $\hat{\Cal{M}}$ over $t$. We denote by $\Cal{M}$ the family of all pairs $(t,[\Sigma])$ where $[\Sigma]$ denotes the equivalence class of $\Sigma$. Likewise, for all $t\in [0,1]$, we denote by $\Cal{M}_t$ the fibre of $\Cal{M}$ over $t$.
\medskip
\noindent For all $t$, we interpret $\Cal{M}_t$ as a smooth Banach manifold (strictly speaking, every relatively compact open subset is an intersection of an infinite family of nested Banach manifolds). We now briefly review the theory of Banach manifolds (see \cite{Klingenberg} for a more detailed description in the $1$ dimensional case). Let $[\Sigma]$ be an element in $\Cal{M}_t$. Let $V_{\Sigma}\subseteq\Cal{M}_t$ be the set of those immersed hypersurfaces which are graphs over $\Sigma$. This is an open subset of $\Cal{M}_t$, which we identify with an open subset $U_{\Sigma}$ of $C_0^\infty(S)$. Let $\Phi_{\Sigma}:U_\Sigma\rightarrow V_\Sigma$ be the canonical identification. $(U_\Sigma,V_\Sigma,\Phi_\Sigma)$ consitutes a smooth chart of $\Cal{M}_t$ which we call the {\bf graph neighbourhood} of $\Sigma$.
\medskip
\noindent We likewise interpret $\Cal{M}$ also as a smooth Banach manifold. As before, let $(t,[\Sigma])$ be an element of $\Cal{M}$, where $\Sigma=(i,(S,\partial S))$. We extend $i$ to a smooth family $(i_s)_{s\in]t-\epsilon,t+\epsilon[}$ such that, for all $s$, $(i_s,\partial S)=\Gamma_s$. Thus, if, for all $s$, we define $\Sigma_s$ by $\Sigma_s = (i_s,(S,\partial S))$, then $(s,[\Sigma_s])_{s\in]t-\epsilon,t+\epsilon[}$ is a smooth family in $\Cal{M}$. Let $V_\Sigma\subseteq\Cal{M}$ be the set of pairs $(s,[\Sigma'])$ where $\Sigma'$ is a graph over $\Sigma_s$. $V_\Sigma$ is an open subset of $\Cal{M}$ which we identify with an open subset, $U_\Sigma$, of $]t-\epsilon,t+\epsilon[\times C_0^\infty(S)$. Let $\Phi_\Sigma:U_\Sigma\rightarrow V_\Sigma$ be the canonical identification. $(U_\Sigma,V_\Sigma,\Phi_\Sigma)$ consitutes a smooth chart of $\Cal{M}$ which we likewise call the {\bf graph neighbourhood} of $\Sigma$. Trivially, this does not depend canonically on $\Sigma$, but also on the choice of smooth family extending $\Sigma$.
\medskip
\noindent Let $(t,\Sigma)$ be an element of $\hat{\Cal{M}}$, where $\Sigma = (i,(S,\partial S))$. The group of smooth diffeomorphisms of $(S,\partial S)$ acts linearly on $C^\infty(S)$. $C^\infty(S)$ therefore defines a bundle $\Cal{E}$ over $\Cal{M}$, whose fibre at $(t,[\Sigma])$ is $C^\infty(S)$. Since the constant functions over $S$ are preserved by the diffeomorphisms of $(S,\partial S)$, these generate a subbundle of $\Cal{E}$ which we identify with $\Cal{M}\times \Bbb{R}$. Likewise, if $(\phi_t)_{t\in[0,1]}\in C^\infty(M)$ is a smooth family of smooth functions, then it defines a section of $\Cal{E}$, which we also denote by $\phi$, given by:
$$
\phi(t,[\Sigma]) = [\phi_t\circ i].
$$
\noindent For all $t$, $\Cal{E}$ restricts canonically to a bundle over $\Cal{M}_t$, which we denote by $\Cal{E}_t$. Let $(U_\Sigma,V_\Sigma,\Phi_\Sigma)$ be a graph neighbourhood of $\Cal{M}_t$ about $\Sigma$. Trivially:
$$
\Cal{E}|_{V_\Sigma} = U_\Sigma\times C^\infty(S).
$$
\noindent This yields a canonical splitting of $T\Cal{E}_t$ over the fibre over $\Sigma$. Since every point in $\Cal{M}_t$ has a canonical graph neighbourhood, we thus obtain a canonical splitting of $T\Cal{E}_t$ which in turn generates a covariant derivative of $\Cal{E}_t$. More explicitely, for every $\Sigma'=(i',(S',\partial S'))\in V_\Sigma$, let $\pi_{\Sigma'}:S'\rightarrow S$ be the canonical projection. A section, $f$, of $\Cal{E}_t$ is covariant constant at $\Sigma$ if and only if there exists a function $f_0\in C^\infty(S)$ such that, up to second order around $\Sigma$:
$$
f_{\Sigma'} = f_0\circ\pi_{\Sigma'}.
$$
\noindent We advise the reader unfamiliar with the theory of Banach manifolds not to trouble himself with the details of this construction. In the sequel, it suffices to know that, locally, $\Cal{E}_t$ behaves like the constant bundle $U_\Sigma\times C^\infty(S)$ and it is not really necessary to have an explicit choice of splitting of $\Cal{E}$.
\medskip
\noindent We define the {\bf Gauss curvature mapping}, $K$, to be the mapping that associates to every element $(t,[\Sigma])$, where $\Sigma = (i,(S,\partial S))$, the function $f\in C^\infty(S)$ whose value at the point $p\in S$ is the Gaussian curvature of $\Sigma$ at $p$. $K$ defines a smooth section of $\Cal{E}$ over $\Cal{M}$.
\medskip
\noindent We determine a formula for the covariant derivative, $\nabla K$, of $K$ with respect to the canonical splitting of $\Cal{E}_t$. Let $\Sigma = (i,(S,\partial S))$ be an element of $\Cal{M}_t$. Let $\msf{N}$ be the outward pointing unit normal vector field over $\Sigma$. Let $R$ be the Riemann curvature tensor of $M$. We define the operator $W$ acting on sections of $TS$ by:
$$
W\cdot X = R_{\msf{N}X}\msf{N}.
$$
\proclaim{Lemma \nextprocno}
\noindent With respect to the canonical splitting, identifying $T_{[\Sigma]}\Cal{M}_t$ with $C_0^\infty(S)$:
$$
\nabla_f K = K\opTr(A^{-1}(W-A^2))f - K\opTr(A^{-1}\opHess(f)),
$$
\noindent where $A$ is the shape operator of $\Sigma$.
\endproclaim
\proof See Proposition $3.1.1$ of \cite{LabI}.\qed
\medskip
\noindent This yields the following result, which will be of use in the sequel:
\proclaim{Corollary \nextprocno}
\noindent $\nabla K$ is a second order linear differential operator. Moreover:
\medskip
\myitem{(i)} if $\Sigma$ is strictly convex, then $\nabla K$ is elliptic; and
\medskip
\myitem{(ii)} when $\opTr(A^{-1}(W-A^2))>0$, $\nabla K$ has trivial kernel.
\endproclaim
\proclabel{CorInvertibility}
\remark In particular, if the sectional curvature of $M$ is bounded above by $-1$ and if $A\leqslant\opId$, then $W-A^2\geqslant 0$ and so, by $(ii)$, $\nabla K$ is invertible.
\medskip
\proof $(i)$ is immediate. $(ii)$ follows by the Maximum Principal.\qed
\goodbreak
\newhead{Locally Convex Hypersurfaces}
\noindent Let $M^{n+1}$ be a Riemannian manifold. A {\bf locally convex hypersurface} in $M$ is a pair $\Sigma=(i,S^n)$ where $S$ is an $n$-dimensional topological manifold and $i:S\rightarrow M$ is a continuous map such that, for all $p\in S$, there exists a neighbourhood, $U$, of $p$ in $S$, a convex subset $K\subseteq M$ with non-trivial interior, and an open subset $V\subseteq\partial K$ such that $i$ restricts to a homeomorphism from $U$ to $V$. We refer to such a triplet $(U,V,K)$ as a {\bf convex chart} of $\Sigma$. Pulling back the metric on $M$ through $i$ yields a natural length metric on $\Sigma$ which we denote by $d_\Sigma$. Let $(\Sigma_n)_\ninn = (i_n,S_n)_\ninn$ and $S_0 = (i_0,S_0)$ be convex immersions. We say that $(\Sigma_n)_\ninn$ {\bf converges} to $\Sigma_0$ if and only if:
\headlabel{HeadConvexImmersions}
\medskip
\myitem{(i)} $(S_n,d_{\Sigma_n})_\ninn$ converges to $(S_0,d_{\Sigma_0})$ in the Gromov-Hausdorff sense; and
\medskip
\myitem{(ii)} $(i_n)_\ninn$ converges to $i_0$ locally uniformly.
\medskip
\noindent Let $\Sigma=(i,S)$ and $\Sigma'=(i',S')$ be two locally convex hypersurfaces in $M$. We say that $\Sigma$ and $\Sigma'$ are equivalent if and only if there exists a homeomorphism $\phi:S\rightarrow S'$ such that:
$$
i = i'\circ\phi.
$$
\example Let $K\subseteq M$ be a convex subset with non trivial interior. Then any open subset of $\partial K$ is a locally convex hypersurface.\qed
\medskip
\example Let $\Sigma$ be a smooth hypersurface on $M$. $\Sigma$ is a locally convex hypersurface if and only if its second fundamental form is everywhere non-negative definite.\qed
\medskip
\noindent Suppose now that $M$ is a Hadamard manifold. Let $K\subseteq M$ be a convex set with non-trivial interior. Let $K^o$ be the interior of $K$. We define $\pi_K:M\setminus K^o\rightarrow \partial K$ to be projection onto the closest point in $\partial K$. Let $V\subseteq\partial K$. We call the set $\pi_K^{-1}(V)$ the {\bf end} of $V$, and we denote it by $\Cal{E}(V)$. Trivially, $\Cal{E}(V)$ is foliated by half geodesics leaving points in $V$ in directions normal to $K$. Let $\Sigma$ be a locally convex hypersurface. Let $(U,V,K)$ and $(U',V',K')$ be convex charts of $\Sigma$. Trivially:
$$
\pi_K^{-1}(i(U\minter U')) = \pi_{K'}^{-1}(i(U\minter U')).
$$
\noindent We thus define the {\bf end} of $\Sigma$ to be the manifold (with non-smooth, concave boundary) whose coordinate charts are the ends of the convex charts of $\Sigma$. We denote this manifold by $\Cal{E}(\Sigma)$. $\Cal{E}(\Sigma)$ has the following properties:
\medskip
\myitem{(i)} $\Sigma$ naturally embeds as the boundary of $\Cal{E}(\Sigma)$;
\medskip
\myitem{(ii)} in the complement of $\Sigma$, $\Cal{E}(\Sigma)$ has the structure of a smooth Riemannian manifold with non-positive curvature;
\medskip
\myitem{(iii)} $\Cal{E}(\Sigma)$ is foliated by half geodesics leaving points in $\Sigma$ in directions normal to $\Sigma$; and
\medskip
\myitem{(iv)} there exists a natural embedding $I:\Cal{E}(\Sigma)\rightarrow M$ which restricts to $i$ over $\Sigma$ and which is a local diffeomorphism over the complement of $\Sigma$.
\medskip
\noindent Let $K\subseteq\Cal{E}(\Sigma)$ be a subset of the end of $\Sigma$. Suppose moreover that $K$ contains $\Sigma$ and that $K$ coincides with $\Sigma$ outside a compact set. Let $p$ be a point in $\Cal{E}(\Sigma)\setminus\Sigma$ lying on the boundary of $K$. We say that $K$ is {\bf boundary convex} at $p$ if and only there exists a neighbourhood, $U$, of $p$ in $\Cal{E}(\Sigma)$, a convex subset $K'\subseteq M$ with non trivial interior, and a neighbourhood $V$ of $I(p)$ in $M$ such that $I$ restricts to a homeomorphism from $U$ to $V$, and:
$$
I(K\minter U) = K'\minter V.
$$
\noindent Bearing in mind that, near any point $p\in\Sigma$, $\Cal{E}(\Sigma)$ may always be extended over an open set containing $p$, we extend this definition to also include boundary points lying in $\Sigma$. We then say that $K$ is {\bf boundary convex} if and only if it is boundary convex at $p$ for every $p\in\partial K$. Importantly, the image under $I$ of the boundary of a boundary convex set is a locally convex hypersurface.
\medskip
\noindent We say that a subset $K\subseteq\Cal{E}(\Sigma)$ is {\bf semi-convex} if and only if for every geodesic segment $\gamma:[0,1]\rightarrow\Cal{E}(\Sigma)$ contained within $\Cal{E}(\Sigma)$, if $\gamma(0),\gamma(1)\in K$, then the whole of $\gamma$ is contained in $K$.
\proclaim{Proposition \nextprocno}
\noindent Let $K$ be a subset of the end of $\Sigma$ which contains $\Sigma$ and coincides with $\Sigma$ outside a convex set. If $K$ is semi-convex, then $K$ is boundary convex.
\endproclaim
\proclabel{LemmaBoundaryConvexityIsSemiConvexity}
\proof Let $p\in\partial K$. If $p$ lies in the interior of $\Cal{E}(\Sigma)$, then $K$ is trivially boundary convex at $p$. Suppose therefore that $p\in\Sigma$. Let $(U,V,K')$ be a convex chart of $\Sigma$ at $p$. Let $r>0$ be such that $B_r(p)\subseteq\Cal{E}(U)$. Consider $X=(K'\minter B_r(p))\munion (K\minter B_r(p))$. Let $\gamma:[0,1]\rightarrow B_r(p)$ be a geodesic segment with endpoints in $X$. Let $\gamma'$ be a maximal subsegment of $\gamma$ lying outside $K'\minter B_r(p)$. Since $\Sigma\subseteq K$, the endpoints of $\gamma'$ are contained in $K\minter B_r(p)$. Thus, by semi-convexity, $\gamma'$ is contained in $K$, and therefore also in $X$. It follows that the whole of $\gamma$ is contained in $X$. Since $\gamma$ is arbitrary, $X$ is convex and $K$ is therefore boundary convex at $p$. This completes the proof.\qed
\medskip
\noindent Let $K$ be a semi-convex subset of the end of $\Sigma$ which contains $\Sigma$ and coincides with $\Sigma$ outside a convex set. $(\partial K,I|_{\partial K})$ defines a convex immersion in $M$ which, by abuse of notation, we simply denote by $\partial K$. Let $\Sigma$ and $\Sigma'$ be two locally convex hypersurfaces in $M$. We say that $\Sigma$ is {\bf bounded by} $\hat{\Sigma}'$ (and $\Sigma'$ {\bf bounds} $\Sigma$) if and only if there exists a semi-convex subset, $K\subseteq\Cal{E}(\Sigma)$, which contains $\Sigma$ and which coincides with $\Sigma$ outside a convex set such that $\Sigma'$ is equivalent to $\partial K$. In this case, we often identify $\Sigma'$ with $\partial K$ and thus view it as a subset of $\Cal{E}(\Sigma)$.
\medskip
\example Let $K,K'\subseteq M$ be two convex sets. Then $\partial K$ is bounded by $\partial K'$ if and only if $K\subseteq K'$.\qed
\medskip
\noindent Let $\Sigma=(i,S)$ be a locally convex hypersurface. For $p\in S$, let $N_p\subseteq UM$ be the set of supporting normals of $\Sigma$ at $S$. We define $N\Sigma$ by:
$$
N\Sigma = \munion_{p\in S}N_p.
$$
\noindent $N\Sigma$ defines a $C^0$ immersed submanifold of $UM$ which we call the {\bf normal} of $\Sigma$.
\medskip
\noindent If $\Sigma'$ bounds $\Sigma$, then there exists an upper semi-continuous function $f:N\Sigma\rightarrow[0,\infty[$ such that $\Sigma'$ is the graph of $f$ over $\Sigma$. Moreover, $f$ vanishes outside a compact set. We call $f$ and $\opSupp(f)$ respectively the {\bf graph function} and {\bf graph support} of $\Sigma'$ with respect to $\Sigma$.
\medskip
\noindent The property of boundedness is preserved by passage to limits:
\goodbreak
\proclaim{Lemma \nextprocno}
\noindent Let $(\Sigma_n)_\ninn,\Sigma_0$ and $(\Sigma'_n)_\ninn,\Sigma'_0$ be convex immersions in $M$. Suppose that, for all $n>0$, $\Sigma'_n$ bounds $\Sigma_n$. For all $n>0$, let $f_n$ and $X_n=\opSupp(f_n)$ be the graph function and graph support respectively of $\Sigma'_n$ with respect to $\Sigma_n$. Suppose that there exists $R>0$ and that, for all $n$, there exists a compact set $X_n'\subseteq\Sigma_n$ such that:
\medskip
\myitem{(i)} $f_n\leqslant R$ for all $n>0$;
\medskip
\myitem{(ii)} for all $n>0$, $X_n\subseteq X_n'$; and
\medskip
\myitem{(iii)} $(X_n')_\ninn$ converges to $X_0'$ in the Hausdorff sense,
\medskip
\noindent then $\Sigma'_0$ also bounds $\Sigma_0$.
\endproclaim
\proclabel{LemmaClosednessOfGraphProperty}
\proof For all $n$, let $K_n\subseteq\Cal{E}(\Sigma_n)$ be the semi-convex subset such that $\partial K_n = \Sigma'_n$. The hypotheses on $(f_n)_\ninn$ and $(X_n)_\ninn$ imply that $(K_n)_\ninn$ is uniformly bounded. By compactness of the family of semi-convex sets, $(K_n)_\ninn$ subconverges in the Hausdorf sense to a semi-convex set $K_0\subseteq\Cal{E}(\Sigma_0)$, say. By Proposition \procref{LemmaBoundaryConvexityIsSemiConvexity}, $K_0$ is boundary convex and so $(I|_{\partial K_0},\partial K_0)$ is a locally convex hypersurface. Moreover $(I|_{\partial K_n},\partial K_n)_\ninn$ converges to $(I|_{\partial K_0},\partial K_0)$ in the sense of convex immersions. Since $(I|_{\partial K_n},\partial K_n)=\Sigma_n'$ for all $n$, and since $(\Sigma_n')_\ninn$ converges to $\Sigma_0'$ in the sense of convex immersions, $(I|_{\partial K_0},\partial K_0)$ is equivalent to $\Sigma_0'$. $\Sigma_0'$ therefore bounds $\Sigma_0$, and this completes the proof.\qed
\medskip
\noindent In the sequel, we require a slight variation of this definition. Let $\Sigma = (i,(S,\partial S))$ and $\Sigma' = (i',(S',\partial S'))$ be (smooth) immersed hypersurfaces which are also convex. Let $\msf{N}_\Sigma$ and $\msf{N}_{\Sigma'}$ be the outward pointing normal vector fields over $\Sigma$ and $\Sigma'$ respectively. Let $\msf{N}_{\partial\Sigma}$ be the normal vector field over $\partial\Sigma$ which is tangent to $\Sigma$ and points outwards from $\partial\Sigma$.
\medskip
\noindent Suppose that $\partial\Sigma'=\partial\Sigma=:\Gamma$. We suppose moreover that $\Sigma'$ lies ``locally strictly above'' $\Sigma$ along $\Gamma$: i.e. for all $p\in\Gamma$:
$$
\langle \msf{N}_{\Sigma'},\msf{N}_{\partial\Sigma}\rangle > 0.
$$
\noindent Since $\Sigma'$ is smooth, it may be extended to a (smooth) convex, immersed hypersurface $\tilde{\Sigma}'$ strictly containing $\partial\Sigma'$ in its interior. Let $\Sigma'_c$ denote the collar region of $\tilde{\Sigma}'$ lying outside $\Sigma'$. We define the piecewise smooth immersed hypersurface $\tilde{\Sigma}$ by:
$$
\tilde{\Sigma} = \Sigma\munion\Sigma'_c.
$$
\noindent Since $\tilde{\Sigma}'$ lies locally strictly above $\Sigma$ along $\Gamma$, $\tilde{\Sigma}$ is also a locally convex hypersurface. We now say that $\Sigma'$ {\bf bounds} $\Sigma$ if and only if $\tilde{\Sigma}'$ bounds $\tilde{\Sigma}$.
\medskip
\noindent Suppose that $\Sigma'$ lies locally strictly above $\Sigma$ along $\partial\Sigma$ and bounds $\Sigma$. Let $f$ be the graph function of $\tilde{\Sigma}'$ with respect to $\tilde{\Sigma}$. Let $\pi:\tilde{\Sigma}'\rightarrow N\tilde{\Sigma}$ be the canonical projection. We say that $\Sigma'$ {\bf strictly bounds} $\Sigma$ if and only for all $p\in S'\setminus\partial S'$:
$$
f\circ\pi(p) > 0.
$$
\noindent In this case, the property of strict containment is preserved by small deformations:
\proclaim{Lemma \nextprocno}
\noindent Let $(\Sigma_n)_\ninn,\Sigma_0$ and $(\Sigma'_n)_\ninn,\Sigma'_0$ be smooth, convex, immersed hypersurfaces. Suppose that $\Sigma'_0$ lies locally strictly above $\Sigma_0$ along $\partial\Sigma_0$ and strictly bounds $\Sigma_0$. Suppose moreover that, for all $n$, $\partial\Sigma_n=\partial\Sigma_n'$ and that $(\Sigma_n)_\ninn$ and $(\Sigma'_n)_\ninn$ converge to $\Sigma_0$ and $\Sigma'_0$ respectively. Then, for sufficiently large $n$, $\Sigma'_n$ lies locally strictly above $\Sigma_n$ along $\partial\Sigma_n$ and bounds $\Sigma_n$.
\endproclaim
\proclabel{LemmaOpennessOfGraphProperty}
\proof For all $n$, let $\Sigma_n'=(i_n',(S_n',\partial S_n'))$. For all $n$, $\Cal{E}(\tilde{\Sigma}_n)$ may be extended beyond $\tilde{\Sigma}_n$ to contain a neighbourhood of $\tilde{\Sigma}_n$. Let $\Cal{E}_\opext(\tilde{\Sigma}_n)$ denote this extension. For sufficiently large $N$, $\Sigma'_n$ is contained in $\Cal{E}_\opext(\tilde{\Sigma}_n)$. Let $d_n:S_n'\rightarrow\Bbb{R}$ be the signed distance in $\Cal{E}_\opext(\tilde{\Sigma}_n)$ to $\tilde{\Sigma}_n$. For sufficently large $n$, $d_n$ is smooth, and $(d_n)_\ninn$ converges to $d_0$ in the $C^\infty$ sense. However, $d_0>0$ and $\nabla d\neq 0$ along $\partial\Sigma_0$. Thus, for sufficiently large $n$, $d_n>0$ and so $\Sigma'_n\subseteq\Cal{E}(\tilde{\Sigma}_n)$. This completes the proof.\qed
\goodbreak
\newhead{Convexity in Higher Codimension}
\noindent Let $M^{n+1}$ be a Riemannian manifold. Let $\Gamma^k=(i,(G^k,\partial G^k))\subseteq M$ be a $k$-dimensional immersed submanifold. Let $N\Gamma\subseteq i^*(UM)$ be the bundle of unit normal vectors over $\Gamma$. $N\Gamma$ has spherical fibres of dimension $(n-k)$. For all $\msf{N}_p\in N\Gamma$, let $A_\Gamma(\msf{N}_p)$ be the shape operator of $\Gamma$ with respect to $\msf{N}_p$. In other words, for all vector fields $X$ and $Y$ tangent to $\Gamma$:
\headlabel{HeadConvexityInHigherCodimension}
$$
A_\Gamma(\msf{N}_p)(X,Y) = - \langle\nabla_X Y,\msf{N}_p\rangle.
$$
\noindent For all $p\in\Gamma$, we define define $X_p\subseteq T_p\Gamma$ by:
$$
X_p = \left\{\msf{N}_p\text{ s.t. }A_\Gamma(\msf{N}_p)> 0\right\},
$$
\noindent where, for a matrix, $M$, we write $M>0$ if and only if it is positive definite. Since the set of positive definite matrices is an open convex cone, $X_p$ is a convex subset of $N_p\Gamma$. In particular, it is contained within a hemisphere. We say that $\Gamma$ is {\bf locally strictly convex} at $p$ if and only if $X_p$ is non-empty. We say that $\Gamma$ is {\bf locally strictly convex} if and only if it is locally strictly convex at every point $p\in\Gamma$.
\medskip
\noindent We now consider the case where $\Gamma$ is of codimension $2$, in which case $N\Gamma$ is a circle bundle over $\Gamma$ and, for all $p\in\Gamma$, $X_p$ is an open interval of length at most $\pi$. We define a {\bf convexity orientation} of $\Gamma$ to be a continuous section, $\msf{N}^-$, of $N\Gamma$ over $\Gamma$ such that, for all $p\in\Gamma$:
$$
\msf{N}^-(p) \in \partial X_p.
$$
\noindent We say that $\Gamma$ carries a convexity orientation when such a section exists. A convexity orientation defines an order over $X_p$ in the following manner: we say that, given two vectors, $V_p,V'_p\in X_p$, $V_p$ lies below $V'_p$ if and only if it lies between $\msf{N}^-(p)$ and $V'_p$. Given a convexity orientation, $\msf{N}^-$, we define the section $\msf{N}^+$ such that, for all $p$:
$$
\partial X_p = \left\{\msf{N}^+,\msf{N}^-\right\}.
$$
\noindent We call this vector field the {\bf convexity coorientation} of $\Gamma$.
\medskip
\example If $(\hat{\Sigma},\partial\hat{\Sigma})$ is a strictly convex immersed hypersurface in $M$, then $\Gamma:=\partial\hat{\Sigma}$ is a locally strictly convex, codimension $2$, immersed submanifold. Moreover, $\Gamma$ inherits a convexity orientation from $\hat{\Sigma}$ in the following manner: For $p\in\Gamma$, we identify each unit vector in $N_p\Gamma$ with the (oriented) hyperplane in $T_pM$ normal to that vector. $T_p\hat{\Sigma}$ defines a half-hyperplane with upward pointing unit normal in $X_p$. Let $H_p$ be another (oriented) hyperplane in $X_p$ that is close to $T_p\Sigma$. We say that $H_p$ lies above (resp. below) $T_p\hat{\Sigma}$ if and only if it is a graph over (resp. beneath) $T_p\hat{\Sigma}$. We extend this to an order on $X_p$, and define $\msf{N}^-(p)$ to be the end point of $X_p$ lying below $T_p\hat{\Sigma}$.
\medskip
\noindent More formally, for $p\in\Gamma$, let $E=T_pM/T_p\Gamma$. $E$ is a two dimensional vector space. Moreover, $N_p\Gamma$ projects down to a circle, $S_p$, in $E$. We consider $X_p$ as a subinterval of $S_p$. Let $\msf{N}_p\in X_p$ be the outward pointing exterior normal to $\hat{\Sigma}$ at $p$. $T_p\hat{\Sigma}$ defines a half-hyperplane which projects down to a half line in $E$. This half-line is parallel to the tangent line to $X_p$ at $\msf{N}_p$, and thus defines an orientation on $S_p$ at $\msf{N}_p$. $\msf{N}^-(p)$ is then the boundary point of $S_p$ towards which $T_p\hat{\Sigma}$ points.\qed
\medskip
\noindent Suppose that $\Gamma$ is locally strictly convex with convexity orientation, and suppose that $\partial\Sigma$ is a strictly convex immersed hypersurface such that $\partial\Sigma = \Gamma$. We say that $\Sigma$ is {\bf compatible} with the orientation on $\Gamma$ if and only if the convexity orientation induced on $\Gamma$ by $\Sigma$ coincides with the pre-existing convexity orientation on $\Gamma$.
\goodbreak
\newhead{First Order Upper Bounds}
\noindent Let $M^{n+1}$ be an $(n+1)$-dimensional Riemannian manifold. Let $\Gamma^{n-1}\subseteq M$ be a strictly convex, codimension $2$, embedded submanifold with convexity orientation. Let $\msf{N}^-$ and $\msf{N}^+$ be the convexity orientation and coorientation respectively of $\Gamma$. Let $\Sigma$ be a strictly convex immersed hypersurface in $M$ such that $\partial\Sigma=\Gamma$. Suppose, moreover that $\Sigma$ is compatible with the convexity orientation on $\Gamma$. We denote by $\msf{N}_\Sigma$ the outward pointing unit normal over $\Sigma$.
\headlabel{FirstOrderUpperBounds}
\medskip
\noindent First order bounds near the boundary follow from the following result:
\proclaim{Lemma \nextprocno}
\noindent Choose $\theta>0$. There exists $r>0$, which only depends on $M$, $\Gamma$ and $\theta$ such that if the angle between $\msf{N}_\Sigma$ and $\msf{N}^+$ is always greater than $\theta$, then, for all $p\in\Gamma$, there exists a convex subset $K\subseteq B_r(p)$ such that the connected component of $\Sigma\minter B_r(p)$ containing $p$ is embedded and is a subset of $\partial K$.
\endproclaim
\proclabel{LemmaBoundaryCompactnessOfConvexImmersions}
\proof This follows immediately from Proposition \procref{PropMPEmbeddedHypersurface} (below).\qed
\medskip
\noindent We establish the framework. Choose $p\in\Gamma$. Choose $r_1>0$, and denote the connected component of $\Gamma\minter B_{r_1}(p)$ containing $p$ by $\Gamma_0$. Reducing $r_1$ if necessary, there exists a smooth, embedded, locally strictly convex hypersurface $\hat{\Sigma}\subseteq B_{r_1}(p)$ such that $\partial\hat{\Sigma}\subseteq\partial B_{r_1}(p)$ and $\Gamma\subseteq\hat{\Sigma}$. We may suppose, moreover, that $\hat{\Sigma}$ bounds a convex set, $\hat{K}$, in $B_{r_1}(p)$. In the sequel, we will identify $M$ with $B_{r_1}(p)$, reducing $r_1$ at various stages whenever necessary. We may thus assume that $\Gamma$ divides $\hat{\Sigma}$ into two connected components: $\hat{\Sigma}^+$ and $\hat{\Sigma}^-$ which correspond to the interior and exterior respectively of $\hat{\Sigma}$ with respect to $\Gamma$.
\medskip
\noindent Let $\msf{N}_{\hat{\Sigma}}$ be the unit normal vector field over $\hat{\Sigma}$. We may suppose that $\msf{N}_{\hat{\Sigma}}$ makes an angle of less than $\theta/2$ with $\msf{N}^+_\Gamma$.
\medskip
\noindent Since $\hat{\Sigma}$ is strictly convex, there exists $\epsilon>0$ such that the shape operator of $\hat{\Sigma}$ is bounded below by $\epsilon\opId$. Let $H$ be a strictly convex embedded hypersurface tangent to $\hat{\Sigma}$ at $p$ whose second fundamental form is strictly bounded above by $\delta\opId$, for $\delta<\epsilon/2$. Let $(H_t)_{t\in]-\tau,\tau[}$ be the foliation of $M$ by hypersurfaces equidistant to $H$. We may assume that each leaf of this foliation is embedded, strictly convex and complete with second fundamental form strictly bounded above by $\delta\opId$. Moreover, we may assume that $H_0=H$ meets $\hat{\Sigma}$ at a single point, $p$. Thus, the upward pointing normal of $H_0$ coincides with that of $\hat{\Sigma}$ at this point.
\medskip
\noindent Each leaf of $(H_t)_{t\in]-\tau,\tau[}$ divides $M$ into two connected components, one of which we say lies above the leaf, and the other of which we say lies below the leaf. Recalling section \headref{HeadConvexImmersions}, we say that a subset $K$ of $M$ is {\bf semi-convex} with respect to a leaf $H_t$ if and only if:
\medskip
\myitem{(i)} it lies above that leaf; and
\medskip
\myitem{(ii)} if $\gamma$ is a geodesic segment lying above $H_t$ whose endpoints are elements of $K$, then the whole of $\gamma$ is contained in $K$.
\medskip
\remark Importantly, in contrast to the situation considered in Section \headref{HeadConvexImmersions}, $H_t$ is not contained in $K$. Semi-convexity is therefore no longer necessarily preserved by taking limits. This is a delicate point which will be discussed presently.
\medskip
\noindent We extend $\Sigma$ to a (piecewise smooth) convex immersed hypersurface by adjoining to it $\hat{\Sigma}^-$ and denote the resulting immersed hypersurface by $\tilde{\Sigma}$. For all $t$, let $\tilde{\Sigma}_t$ be the connected component of $\tilde{\Sigma}$ lying above $H_t$ and containing $p$.
\medskip
\noindent Lemma \procref{LemmaBoundaryCompactnessOfConvexImmersions} follows immediately from the following proposition by taking intersections with a small ball about $p$:
\proclaim{Proposition \nextprocno}
\noindent There exists $t_0<0$ (which only depends on $M$, $\Gamma$, $\hat{\Sigma}$, $\theta$ and $r_1$) such that $\tilde{\Sigma}_{t_0}$ is embedded and (along with $H_{t_0}$) bounds a semi-convex set.
\endproclaim
\proclabel{PropMPEmbeddedHypersurface}
\proof This follows immediately from Proposition \procref{PropMPUpperBoundOfInfimum} (below).\qed
\medskip
\noindent Let $T$ denote the set of all $t<0$ such that, for all $s\in]-t,0[$:
\medskip
\myitem{(i)} $\tilde{\Sigma}_s$ is embedded;
\medskip
\myitem{(ii)} $\tilde{\Sigma}_s\subseteq\hat{K}$;
\medskip
\myitem{(iii)} $\tilde{\Sigma}_s$ bounds a semi-convex set above $H_s$; and
\medskip
\myitem{(iv)} $\tilde{\Sigma}_s$ intersects $H_s$ transversally along $\partial\tilde{\Sigma}_s$.
\proclaim{Proposition \nextprocno}
\noindent $T$ is non-empty.
\endproclaim
\proof Since $\tilde{\Sigma}$ is a piecewise smooth, strictly convex immersion, there exists $0<r_2<r_1$ (which {\sl does\/} depend on $\Sigma$) such that the connected component of the intersection of $\tilde{\Sigma}$ with $B_{r_2}(p)$ containing $p$ is embedded and bounds a convex set. The portion of this convex set lying above $H_t$ for $t$ small is trivially semi-convex, and $(i)$, $(ii)$ and $(iii)$ are therefore satisfied for all small $t$ less than $0$, likewise so is $(iv)$, and the result follows.\qed
\medskip
\noindent Let $t_0$ be the infimum of $T$. We will obtain upper bounds for $t_0$, from which Proposition \procref{PropMPEmbeddedHypersurface} will follow. The first step involves proving that $\tilde{\Sigma}_{t_0}$ is transverse to $H_{t_0}$. The main geometric obstacle is the possibility that the outward pointing normal to $\Sigma_{t_0}$ points upwards from $H_{t_0}$. This is dealt with by the following observation:
\proclaim{Proposition \nextprocno}
\noindent For all $\theta'<\theta$, there exists $t_1<0$ (which only depends on $\hat{\Sigma}$, $\Gamma$, $M$, $\theta$ and $\theta'$) such that, if $d$ is the (signed) distance function in $M$ to $H=H_0$, and if $t_0>t_1$, then, throughout $\Sigma_{t_0}$:
$$
\langle\msf{N}_\Sigma,\nabla d\rangle \leqslant \opCos(\theta') < 1.
$$
\endproclaim
\proclabel{PropAngleBddBelow}
\proof Let $\nabla$ and $\nabla^\Sigma$ denote the Levi-Civita covariant derivatives over $M$ and $\Sigma$ respectively. Define the function $\phi:\Sigma\rightarrow\Bbb{R}$ by:
$$
\phi = \langle\msf{N},\nabla d\rangle.
$$
\noindent Let $A$ be the shape operator of $\Sigma$. If $X$ is a vector field over $\Sigma$, then:
$$\matrix
X\phi \hfill&= \langle\nabla_X\msf{N}_\Sigma,\nabla d\rangle + \langle\msf{N}_\Sigma,\nabla_X\nabla d\rangle\hfill\cr
&= \langle A\cdot X,\nabla d\rangle + \opHess(d)(\msf{N}_\Sigma,X)\hfill\cr
&= \langle A\cdot X,\nabla^\Sigma d\rangle + \opHess(d)(\msf{N}_\Sigma,X).\hfill\cr
\endmatrix$$
\noindent The final line follows since the normal component of $A\cdot X$ vanishes. Now let $X=\nabla^\Sigma d/\|\nabla^\Sigma d\|^2$. Since $A$ is positive definite:
$$
X\phi \geqslant \opHess(d)(\msf{N}_\Sigma,X).
$$
\noindent Since $d$ is the distance to a hypersurface, $\|\nabla d\|=1$ is constant, and so $\opHess(d)(\nabla d,\cdot)$ vanishes. Thus, if $\msf{N}_0$ denotes the component of $\msf{N}_\Sigma$ tangent to the foliation, $(H_t)_{t\in]-\tau,\tau[}$, of level subsets of $d$, then:
$$
X\phi \geqslant \opHess(d)(\msf{N}_0,X).
$$
\noindent However:
$$
\|\msf{N}_0\|^2 = 1 - \langle\msf{N}_\Sigma,\nabla d\rangle^2 = \|\nabla^\Sigma d\|^2.
$$
\noindent Thus $\msf{N}_0/\|\nabla^\Sigma d\|$ has norm equal to $1$. Since the shape operator of $H_t$ is bounded above by $\delta\opId$ for all $t$, the norm of $\opHess(d)$ is also bounded above by $\delta$. Thus:
$$
X\phi \geqslant -\delta.
$$
\noindent However:
$$
Xd = \langle X,\nabla^\Sigma d\rangle = 1.
$$
\noindent Thus, if $\gamma:[0,\tau]\rightarrow\Sigma$ is an integral curve of $X$ starting at $q$, then $\gamma(s)$ meets $\Gamma$ for some $s\leqslant\left|t_0\right|$. There therefore exists $q'\in\Gamma$ such that:
$$\matrix
&\phi(q') \hfill&\geqslant \phi(q) - \epsilon\left|t_0\right|/2 \hfill\cr
\Rightarrow\hfill&\phi(q) \hfill&\leqslant \phi(q') + \epsilon\left|t_0\right|/2,\hfill\cr
& &\leqslant \phi(q')+ \epsilon\left|t_1\right|/2.\hfill\cr
\endmatrix$$
\noindent Choosing $t_1$ sufficiently small, for all $q'\in\Gamma_{t_1}$:
$$
\phi(q') + \epsilon\left|t_1\right|/2\leqslant\opCos(\theta').
$$
\noindent The result follows.\qed
\proclaim{Proposition \nextprocno}
\noindent There exists $t_1<0$ (which only depends on $M$, $\hat{\Sigma}$, $\theta$ and $r_1$) such that, if $t_0>t_1$, then $\tilde{\Sigma}_{t_0}$ intersects $H_{t_0}$ transversally along $\partial\tilde{\Sigma}_{t_0}$.
\endproclaim
\proclabel{PropSecondStep}
\proof Suppose the contrary. Choose $q\in\partial\tilde{\Sigma}_{t_0}$ such that $\tilde{\Sigma}$ is tangent to $H_{t_0}$ at $q$. The normal to $\tilde{\Sigma}$ at $q$ either points downwards into $H_{t_0}$ or upwards from $H_{t_0}$. By reducing $r_1$ if necessary, we may assume that the normal does not point upwards over $\hat{\Sigma}^-\setminus\left\{p\right\}$. By Proposition \procref{PropAngleBddBelow}, for $t_1$ sufficiently small, the normal doesn't point upwards over $\Sigma$ either, and it therefore does not point upwards anywhere over $\tilde{\Sigma}_{t_0}$.
\medskip
\noindent We now show that the normal cannot point downwards. Since $\tilde{\Sigma}$ and $H_{t_0}$ are strictly convex with opposing normals, they meet at a single point. For $t>t_0$, let $\tilde{\Sigma}'_t$ denote the connected component of $\tilde{\Sigma}$ lying below $H_t$ containing $q$. Since $\tilde{\Sigma}$ is piecewise smooth, for $t$ sufficiently close to $t_0$, $\tilde{\Sigma}'_t$ is topologically a ball whose boundary is an embedded topological sphere in $H_t$. Moreover, $\partial\tilde{\Sigma}'_t$ is a subset of $\partial\tilde{\Sigma}_t$. However, for all $s\in]t_0,0[$, $\tilde{\Sigma}_s$ is transverse to $H_s$ and does not self intersect. $\partial\tilde{\Sigma}_t$ is thus also an embedded topological sphere. It follows that $\partial\tilde{\Sigma}_t$ and $\partial\tilde{\Sigma}'_t$ coincide, and $\tilde{\Sigma}$ is therefore an embedded topological sphere lying above $H_{t_0}$. $\Gamma$ is therefore not contained in $\tilde{\Sigma}$, which is absurd, and thus the normal to $\tilde{\Sigma}'_{t_0}$ does not point downwards, and this completes the proof.\qed
\medskip
\noindent The next step uses the fact that $K_{t_0}$, being the limit of a sequence of semi-convex sets, is also semi-convex. Despite being an intuitive result, its proof is rather technical, and is deferred to Section \headref{HeadSemiConvexity}.
\proclaim{Proposition \nextprocno}
\noindent There exists $t_1<0$ (which only depends on $M$, $\hat{\Sigma}$, $\theta$ and $t_1$), such that, if $t_0>t_1$, then $\partial\tilde{\Sigma}_{t_0}$ is embedded in $H_{t_0}$ and bounds an open set.
\endproclaim
\proclabel{PropMPNoSelfIntersection}
\proof For $t>t_0$, let $K_t$ be the semi-convex set bounded by $\tilde{\Sigma}_t$ and $N_t$. By Proposition \procref{PropClosednessOfSemiConvexity}, $K_{t_0}$ is also semi-convex. By the preceeding proposition, $\tilde{\Sigma}_{t_0}$ is transverse to $H_{t_0}$ along $\partial\tilde{\Sigma}_{t_0}$. It follows that $\partial\tilde{\Sigma}_{t_0}$ is a (piecewise smooth) immersed submanifold of $H_{t_0}$. Suppose it is not embedded. Since $\partial\Sigma_t$ is embedded for all $t>t_0$, there exist two open subsets $\Sigma'_1,\Sigma'_2\subseteq\Sigma_{t_0}$ such that:
\medskip
\myitem{(i)} $\Sigma'_1$ and $\Sigma'_2$ are embedded; and
\medskip
\myitem{(ii)} $\Sigma'_1\minter H_{t_0}$ and $\Sigma'_2\minter H_{t_0}$ meet tangentially at some point $p$.
\medskip
\noindent Since $\tilde{\Sigma}_{t_0}$ bounds $K_{t_0}$, the hypersurfaces $\Sigma'_1$ and $\Sigma'_2$ divide a neighbourhood of $p$ above $H_{t_0}$ into three (roughly) wedge-shaped open sets. Consider the central one of these three wedges. It is either a subset of $K_{t_0}$ or a subset of its complement. If it is a subset of $K_{t_0}$, then we say that $\Sigma'_1\minter H_{t_0}$ and $\Sigma'_2\minter H_{t_0}$ lie on each others interior. Otherwise, we say that they lie on each others exterior.
\medskip
\noindent Suppose that $\Sigma'_1\minter H_{t_0}$ and $\Sigma'_2\minter H_{t_0}$ lie on each others interior. Let $P_1$ and $P_2$ be the respective tangent hyperplanes of $\Sigma'_1$ and $\Sigma'_2$ at $p$. We identify these with their images under the exponential map. $P_1$ and $P_2$ do not coincide. Indeed, suppose the contrary. By strict convexity, the interiors of $\Sigma'_1$ and $\Sigma'_2$ coincide in a single point. This point is contained in $K_{t_0}$. However, $K_{t_0}$ is connected and also contains $p$, which is absurd and the assertion follows.
\medskip
\noindent By convexity, near $p$, $\Sigma'_1$ lies above $P_1$ and $\Sigma'_2$ lies above $P_2$. However, the region lying above both $P_1$ and $P_2$ forms a wedge making an angle at $p$ strictly greater than $0$ and strictly less than $\pi$. In particular, $\Sigma'_1$ and $\Sigma'_2$ intersect transversally at $p$. They therefore also intersect over a hypersurface contained inside this wedge. However, since $H_{t_0}$ is strictly convex, this wedge lies strictly above $H_{t_0}$, and therefore $\Sigma'_1$ and $\Sigma'_2$ also meet at some point above $H_{t_0}$. This contradicts the hypothesis that $\partial\tilde{\Sigma}_t$ is an embedded submanifold of $H_t$ for all $t>t_0$. It follows that these two submanifolds do not lie on each others interior.
\medskip
\noindent Suppose that $\Sigma'_1\minter H_{t_0}$ and $\Sigma'_2\minter H_{t_0}$ lie on each others exterior. Let $\gamma$ be a geodesic arc, tangent to $H_{t_0}$ at $p$ and normal to the common tangent space of $\Sigma'_1\minter H_{t_0}$ and $\Sigma'_2\minter H_{t_0}$. Near $p$, $\gamma$ lies above $H_{t_0}$ and has endpoints inside $K_{t_0}$. Moving $\gamma$ upwards slightly yields a geodesic arc lying above $H_{t_0}$, having endpoints inside $K_{t_0}$ whilst itself not being contained within $K_{t_0}$. This contradicts semi-convexity. It follows that these two submanifolds do not lie on each others exterior, and this completes the proof.\qed
\proclaim{Proposition \nextprocno}
\noindent There exists $t_1<0$ (which only depends on $M$, $\hat{\Sigma}$, $\theta$ and $r_1$) such that, if $t_0>t_1$, then $\tilde{\Sigma}_t\subseteq\hat{K}$.
\endproclaim
\proclabel{PropFinalStep}
\proof By Proposition \procref{PropMPNoSelfIntersection}, for $t$ sufficiently small, $\Sigma_t$ does not intersect $\hat{\Sigma}_t^-$. By Proposition \procref{PropMPEnvelopeLiesInsideSigma}, semi-convexity and the hypotheses on $\Sigma$ along the boundary, for all sufficiently small $t$, $\Sigma_t$ does not intersect $\hat{\Sigma}_t^+$. It is therefore contained within the set bounded by $\hat{\Sigma}_t$ and $H_t$, and the result follows.\qed
\proclaim{Proposition \nextprocno}
\noindent There exists $t_1<0$ (which only depends on $M$, $\hat{\Sigma}$, $\theta$ and $r_1$) such that, if $t_0>t_1$, then $t_0$ cannot be the infimum of $T$.
\endproclaim
\proclabel{PropMPUpperBoundOfInfimum}
\proof Let $t_1$ be as in Propositions \procref{PropAngleBddBelow}, \procref{PropSecondStep}, \procref{PropMPNoSelfIntersection} and \procref{PropFinalStep} and suppose that $t_0>t_1$. $\partial\tilde{\Sigma}_{t_0}$ is embedded, is transverse to $H_{t_0}$, and is bounded away from $B_{r_1}(p)$. Thus, for all $t<t_0$ sufficiently close to $t_0$, $\tilde{\Sigma}_{t}$ is embedded, is contained in $\hat{K}$, meets $H_t$ transversally, and, along with $H_t$, bounds a subset of $B_{r}(p)$. For all $t$, let $K_t$ be the closure of this subset. It thus remains to show that $K_t$ is semi-convex for all $t$ sufficiently close to $t_0$.
\medskip
\noindent Suppose that there exists a sequence $(t_n)_\ninn<t_0$ converging to $t_0$ such that, for all $n$, $K_{t_n}$ is not semi-convex. Then, for all $n$, there exists $p_n,q_n\in K_n:=K_{t_n}$ and a geodesic arc $\gamma_n$ such that:
\medskip
\myitem{(i)} $p_n$ and $q_n$ are the endpoints of $\gamma_n$;
\medskip
\myitem{(ii)} $\gamma_n$ lies above $H_{t_n}$; and
\medskip
\myitem{(iii)} there exists a point $r_n\in\gamma_n$ which lies outside $K_n$.
\medskip
\noindent Without loss of generality, $(p_n)_\ninn$, $(q_n)_\ninn$, $(\gamma_n)_\ninn$ and $(r_n)_\ninn$ converge to $p_0$, $q_0$, $\gamma_0$ and $r_0$ respectively. Trivially, $\gamma_0$ lies above $H_{t_0}$. Suppose first that $p_0\neq q_0$. Suppose that $r_0$ does not coincide with either of the endpoints. Since $K_{t_0}$ is semi-convex, $\gamma_0$ lies inside $K_{t_0}$. $r_0$ therefore lies on the boundary of $K_{t_0}$, and $\gamma_0$ is therefore an interior tangent to $\tilde{\Sigma}_{t_0}$ at this point, which contradicts local strict convexity. Likewise, if $r_0$ coincides with an end point, $p_0$, say, then $\gamma_0$ is contained inside $K_{t_0}$ and points outwards (or is tangent) to $\tilde{\Sigma}_{t_0}$ at $p_0$, which also contradicts local strict convexity and semi-convexity. It follows that $p_0$ and $q_0$ coincide.
\medskip
\noindent If $p_0=q_0$ is an interior point of $K_{t_0}$, then $\gamma_n$ trivially lies inside $K_{t_0}$ for all sufficiently large $n$. Suppose therefore that $p:=p_0=q_0$ is a boundary point of $\tilde{\Sigma}_{t_0}$. By local strict convexity, there exists a neighbourhood of $\tilde{\Sigma}$ about $p$ which lies on the boundary of a convex set, $X$. For all $n$, the intersection of $X$ with the region lying above $H_{t_n}$ is a subset of $K_n$. However, for sufficiently large $n$, $p_n$ and $q_n$ both lie in $X$. For all such $n$, $\gamma_n$ is contained within $X$ and therefore within $K_n$, which is absurd.
\medskip
\noindent There therefore exists $\epsilon$ such that, for all $t>t_0-\epsilon$, $\tilde{\Sigma}_t$ satisfies the hypotheses defining $T$, and therefore $t\in T$. This is absurd and the result follows.\qed
\goodbreak
\newhead{Parabolic Limits}
\noindent Let $M^{n+1}$ be an $(n+1)$-dimensional Riemannian manifold. Let $\hat{\Sigma}^n$ be a locally strictly convex immersed hypersurface in $M$. Let $\Gamma\subseteq\hat{\Sigma}$ be an embedded hypersurface. Let $\epsilon>0$ be such that the shape operator of $\hat{\Sigma}$ is everywhere bounded below by $\epsilon$. Choose $p\in\Gamma$. Let $H^n$ be a strictly convex embedded hypersurface of $M$ which is an exterior tangent to $\hat{\Sigma}$ at $p$. Let $\delta>0$ be such that the shape operator of $H$ is everywhere bounded above by $\delta$ and suppose that $\delta<\epsilon/2$. For simplicity, we assume throughout the rest of this section that the shape operator of $H$ at $p$ is equal to $\delta\opId$. The general case is similar.
\headlabel{HeadParabolicLimits}
\medskip
\noindent Let $d$ be the signed distance function in $M$ to $H$. In particular, for $q\in\hat{\Sigma}$ near $p$, $d(q)\leqslant 0$. For all $t$, let $H_t$ be the level hypersurface at distance $t$ from $H$. For small $t<0$, let $\hat{\Sigma}_t$ and $\Gamma_t$ be the connected components of $\hat{\Sigma}$ and $\Gamma$ respectively lying above $H_t$ and containing $p$, and let $\hat{K}_t$ denote the compact set bounded by $\hat{\Sigma}_t$ and $H_t$. For small $t$, $\Gamma_t$ divides $\hat{\Sigma}_t$ into two components, which we denote by $\hat{\Sigma}_t^+$ and $\hat{\Sigma}_t^-$.
\medskip
\noindent Choose $t_0<0$. Let $(p_n)_\ninn\in\hat{K}_{t_0}$ be a sequence converging to $p$. We consider a geodesic chart for $H$ about $p$, and thus identify a neighbourhood of $p$ in $H$ with a neighbourhood of $0$ in $T_pH$. Let $(e_1,...,e_n)$ be an orthonormal basis for $T_pH$. There exists $r>0$ such that $\hat{\Sigma}$ is the graph of a function, $f$ over $B_r(p)$. By Taylor's Theorem, with respect to $(e_1,...,e_n)$:
$$
f(x) = -\langle x|A|x\rangle + O(\|x\|^3),
$$
\noindent where $A$ is a positive definite matrix. With respect to these coordinates, for all $n$, $p_n=(q_n,s_n)$, where $q_n\in T_pH$ and $s_n<0$. For all $n$, we define $\hat{f}_n:B_{r/\sqrt{\left|s_n\right|}}(p)\rightarrow]-\infty,0[$ and $\hat{q}_n\in T_pH$ by:
$$
\hat{f}_n(x) = f(\sqrt{\left|s_n\right|}x)/\left|s_n\right|,\qquad \hat{q}_n = q_n/\sqrt{\left|s_n\right|}.
$$
\noindent Trivially, $(\hat{f}_n)_\ninn$ converges in the $C^\infty_\oploc$ sense over $T_pH$ to $\hat{f}_0$, where:
$$
\hat{f}_0(x) = -\langle x|A|x\rangle.
$$
\noindent Moreover, for all $n$, since $p_n\in\hat{K}_{s_n}$:
$$\matrix
&\left|f(q_n)\right| \hfill&\leqslant \left|s_n\right|\hfill&\cr
\Rightarrow\hfill&\mlimsup_{n\rightarrow\infty}\frac{\epsilon}{2}\|\hat{q}_n\|^2\hfill
&\leqslant\mlimsup_{n\rightarrow\infty}\frac{1}{\left|s_n\right|}\left|f(q_n)\right|\hfill&\leqslant 1.\cr
\endmatrix$$
\noindent There thus exists $\hat{q}_0\in H$ towards which $(\hat{q}_n)_\ninn$ subconverges. In particular, $\hat{f}_0(q_0)\geqslant -1$. We call $(\hat{f}_0,\hat{q}_0)$ a {\bf parabolic limit} of $(\hat{\Sigma},q_n)_\ninn$.
\medskip
\noindent Likewise, if we suppose that $(e_1,...,e_{n-1})$ is tangent to $\Gamma$ at $p$, then, reducing $r$ if necessary, the projection of $\Gamma$ onto $H$ is the graph of some function $g$, over the space spanned by $(e_1,...,e_{n-1})$. For all $n$, we define $\hat{g}_n:B_{r/\sqrt{\left|s_n\right|}}(p)\rightarrow\Bbb{R}$ by:
$$
\hat{g}_n(x') = g(\sqrt{\left|s_n\right|}x')/\sqrt{\left|s_n\right|}.
$$
\noindent Trivially, $(\hat{g}_n)_\ninn$ subconverges in the $C^\infty_\oploc$ sense over the space spanned by $(e_1,...,e_{n-1})$ to $\hat{g}_0:=0$. It follows that the parabolic limit of $\Gamma_{s_n}$ is the intersection of the graph of $\hat{f}_0$ with a vertical hyperplane in $\Bbb{R}^n\times\Bbb{R}$.
\medskip
\noindent For $p\in M$, we call a {\bf geodesic hyperplane} at $p$ an immersed hypersurface consisting of geodesics passing through $p$. Explicitely, $P\subseteq M$ is a geodesic hyperplane if and only if there exists a hyperplane $H\subseteq T_pM$ such that:
$$
P=\left\{\opExp(V_p)\text{ s.t. }V_p\in H\right\}.
$$
\noindent For all $n$, let $P_n$ be the geodesic hyperplane tangent to $H_{s_n}$ at $p_n$. Reducing $r$ if necessary, $P_n$ is the graph of the function $\phi_n:B_r(q_n)\rightarrow\Bbb{R}$, where, by convexity:
$$
-s_n\leqslant \phi_n(x)\leqslant -s_n + \langle x-q_n|B_n|x-q_n\rangle + O(\|x\|^3),
$$
\noindent where $(B_n)_\ninn$ converges to $\delta\opId$. For all $n$, we define $\hat{\phi}_n$ by:
$$
\hat{\phi}_n(x) = \phi_n(\sqrt{\left|s_n\right|}x)/\left|s_n\right|.
$$
\noindent $(\hat{\phi}_n)_\ninn$ converges in the $C^\infty_\oploc$ sense over $T_pH$ to $\hat{\phi}_0$ where:
$$
\hat{\phi}_0(x) = \delta\|x-\hat{q}_0\|^2 - 1.
$$
\noindent Thus, the parabolic limit of the geodesic hyperplanes tangent to $H_{s_n}$ at $p_n$ is a paraboloid on $(\hat{q}_0,-1)$. Finally, in like manner, the parabolic limit of a sequence of geodesics tangent to $H_{s_n}$ at $p_n$ is the intersection of this paraboloid with a vertical plane in $\Bbb{R}^n\times\Bbb{R}$.
\medskip
\noindent Parabolic limits are of use in obtaining technical results concerning $\Sigma$.
\proclaim{Proposition \nextprocno}
\noindent For $\delta$ sufficiently small, there exists $t_0<0$ (which only depends on $\hat{\Sigma}$, $\Gamma$, $H$ and $M$) such that, for all $q\in K_{t_0}\setminus\left\{p\right\}$, if $t>t_0$ is such that $q\in H_t$, and if $P$ is the geodesic hyperplane tangent to $H_t$ at $q$, then:
\medskip
\myitem{(i)} $P$ intersects $\hat{\Sigma}_t$ transversally; and
\medskip
\myitem{(ii)} $P$ intersects $\Gamma_t$ transverally.
\endproclaim
\proclabel{CorTransversalityOfFixedPart}
\proof $(i)$ Suppose the contrary. Let $(p_n)_\ninn\in\hat{K}_{t_0}$ be a sequence converging to $p$. For all $n$, let $s_n<0$ be such that $p_n\in H_{s_n}$ and let $P_n$ be the geodesic hyperplane tangent to $H_{s_n}$ at $p_n$. Trivially, $P_n$ intersects $\hat{\Sigma}_n$ non-trivially for all $n$. Suppose that, for all $n$, $P_n$ is tangent to $\hat{\Sigma}_{s_n}$ at some point. It follows that the parabolic limit of $(P_n)_\ninn$ is tangent to the parabolic limit of $(\hat{\Sigma}_{s_n})_\ninn$ at some point. This is absurd, and the first assertion follows.
\medskip
\noindent $(ii)$ Suppose the contrary. Let $(p_n)_\ninn\in\hat{\Sigma}$ be a sequence converging towards $p$. For all $n$, let $s_n<0$ be such that $p_n\in H_{s_n}$, let $\Gamma_n=\Gamma_{s_n}$ and let $P_n$ be the geodesic hyperplane tangent to $H_{s_n}$ at $p_n$. We suppose that, for all $n$:
$$
P_n\minter\Gamma_n =\emptyset.
$$
\noindent The parabolic limit of $P_n$ intersects the parabolic limit of $\Gamma_n$ transversally. Thus, for sufficiently large $n$, $P_n\minter\Gamma_n\neq\emptyset$, which is absurd. It follows that, for $t_0$ sufficiently small, $P$ intersects $\Gamma_{t_0}$. Transversality follows as in the proof of part $(i)$, and this completes the proof.\qed
\proclaim{Proposition \nextprocno}
\noindent Choose $\theta\in]0,\pi/2[$. For $\delta$ sufficiently small, there exists $t_0<0$ (which only depends on $\hat{\Sigma}$, $\Gamma$, $H$, $M$ and $\theta$) such that for $t>t_0$ and for all $q\in\hat{\Sigma}_t\minter H_t$, there exists a geodesic segment, $\gamma$, joining $q$ to $\Gamma$ such that the hyperplane spanned by $\partial_t\gamma$ and $T\Gamma$ at the point of intersection of $\gamma$ with $\Gamma$ makes an angle strictly less than $\theta$ with $T\hat{\Sigma}$.
\endproclaim
\proclabel{PropMPEnvelopeLiesInsideSigma}
\proof Suppose the contrary. Let $(p_n)_\ninn\in\hat{\Sigma}$ be a sequence converging to $p$, and let $\hat{p}_0$ be its parabolic limit. For all $n$, let $s_n<0$ be such that $p_n\in H_{s_n}$ and let $\gamma_n$ be a geodesic segment tangent to $H_{s_n}$ at $p$ and terminating in $\Gamma_{s_n}$. Suppose that, for all $n$, the hyperplane spanned by $\partial_t\gamma_n$ and $T\Gamma_{s_n}$ at the point of intersection of $\gamma_n$ with $\Gamma_{s_n}$ makes an angle of at least $\theta$ with $T\hat{\Sigma}$. Let $\hat{\gamma}_0$ and $\hat{\Gamma}_0$ be the parabolic limits of $(\gamma_n)_\ninn$ and $(\Gamma_{s_n})_\ninn$ respectively. Then, at its point of intersection with $\hat{\Gamma}_0$, $\hat{\gamma}_0$ is tangent to the vertical hyperplane containing $\hat{\Gamma}_0$. $\hat{\gamma}_0$ is thus entirely contained in this vertical hyperplane. It follows that every parabolic limit of every sequence of geodesic segments joining $(p_n)_\ninn$ to $\Gamma$ is contained in the vertical hyperplane containing $\hat{\Gamma}_0$. When $\hat{p}_0\notin\hat{\Gamma}_0$, this is trivially absurd. When $\hat{p}_0\in\hat{\Gamma}_0$, there exists a parabolic limit of such geodesic segments which is normal to the hyperplane containing $\hat{\Gamma}_0$, which is also absurd. The result follows.\qed
\proclaim{Proposition \nextprocno}
\noindent For $\delta$ sufficiently small, there exists $t_0<0$ (which only depends on $\hat{\Sigma}$, $\Gamma$, $H$ and $M$) such that, for $t>t_0$, if $\gamma$ is a geodesic segment lying in $\hat{K}_{t_0}$ such that:
\medskip
\myitem{(i)} $\gamma$ is tangent to $H_t$; and
\medskip
\myitem{(ii)} the endpoints of $\gamma$ both lie in $\Gamma$,
\medskip
\noindent then there exists a sequence of geodesic segments $(\gamma_n)_\ninn$ converging to $\gamma$ such that, for all $n$:
\medskip
\myitem{(i)} $\gamma_n$ is tangent to $H_t$; and
\medskip
\myitem{(ii)} the end points of $\gamma_n$ lie in $\hat{\Sigma}_t^-$.
\endproclaim
\proclabel{PropGeodesicsWithBothEndsInBoundary}
\proof Suppose the contrary. Let $(p_n)_\ninn\in\hat{K}_{t_0}$ be a sequence converging to $p$. For all $n$, let $s_n<0$ be such that $p_n\in H_{s_n}$ and let $\gamma_n$ be a geodesic segment tangent to $H_{s_n}$ at $p_n$ with both end points in $\Gamma_{s_n}$. We suppose that, for all $n$, there exists $\epsilon_n>0$ such that if $q_n\in H_{s_n}$ is such that $d(q_n,p_n)<\epsilon_n$, then no geodesic segment tangent to $H_{s_n}$ at $q_n$ has both endpoints in $\hat{\Sigma}^-_{s_n}$. Let $\hat{\gamma}_0$, $\hat{p}_0$, $\hat{\Gamma}_0$ and $\hat{\Sigma}^-_{s_n}$ be the parabolic limits of $(\gamma_n)_\ninn$, $(p_n)_\ninn$, $(\Gamma_n)_\ninn$ and $(\hat{\Sigma}^-_{s_n})_\ninn$ respectively. Let $\hat{V}_0$ be the horizontal unit vector at $(p_0,-1)$ normal to the vertical hyperplane containing $\hat{\Gamma}_0$ and pointing towards $\hat{\Sigma}_0^-$. For all $n$, let $V_n$ be a unit vector tangent to $H_{s_n}$ at $p_n$ and suppose that $\hat{V}_0$ is the parabolic limit of $(V_n)_\ninn$. For all $n$, let $\eta_n:\Bbb{R}\rightarrow H_{s_n}$ be the geodesic in $H_{s_n}$ such that:
$$
\partial_t\eta_n(0) = V_n,
$$
\noindent and let $X_n$ be the parallel transport of $\partial_t\gamma_n(0)$ along $\eta_n$ (with respect to the Levi-Civita covariant derivative of $H_{s_n}$). Let $\opExp$ be the exponential map of $M$ and for all $n$ define:
$$
\phi_{n,t}(s) = \opExp(sX_n(t)).
$$
\noindent If $\hat{X}_0$ is the unit tangent vector to $\hat{\gamma}_0$ at $\hat{p}_0$, then the parabolic limit of $(\phi_n)_\ninn$ is $\hat{\phi}_{0,t}(s)$, where:
$$
\hat{\phi}_{0,t}(s) = (s\hat{X}_0 + t\hat{V}_0,\delta s^2 - 1).
$$
\noindent The intersection of this family with $\hat{\Sigma}^-_0$ is transverse to $\hat{\Gamma}_0$ at the intersection of $\hat{\gamma}_0$ with $\hat{\Gamma}_0$. Thus, for sufficiently large $n$, and sufficiently small $t$, the two endpoints of the geodesic segment $s\mapsto\phi_{n,t}(s)$ both lie in $\hat{\Sigma}_0^-$. This is absurd, and the result follows.\qed
\goodbreak
\newhead{Semi-Convexity}
\noindent In this section we show that the property of being semi-convex is preserved after taking limits. Using the same notation as in the Section \headref{FirstOrderUpperBounds}, we show:
\headlabel{HeadSemiConvexity}
\proclaim{Proposition \nextprocno}
\noindent There exists $t_1<0$ (which only depends on $M$, $\hat{\Sigma}$, $\theta$ and $r_1$) such that, if $t_0\geqslant t_1$, then $\tilde{\Sigma}_{t_0}$ bounds a semi-convex set above $H_{t_0}$.
\endproclaim
\proclabel{PropClosednessOfSemiConvexity}
\noindent For $p\in M$, if $P$ is a geodesic hyperplane at $p$ (see Section \headref{HeadParabolicLimits}), then we say that two points $q_1,q_2\in P$ are {\bf coaxial} if and only if they both lie on the same radial geodesic on opposite sides of $p$. We require the following technical result:
\proclaim{Lemma \nextprocno}
\noindent Choose $\varphi>0$. Let $K\subseteq M$ be compact. There exists $r>0$ (which only depends on $\varphi$ and $K$) such that, if $P$ is a geodesic hyperplane at $p\in K$, if $q_1,q_2\in P$ are coaxial points and if $X$ is a Jacobi field over the geodesic joining $q_1$ to $q_2$ such that:
\medskip
\myitem{(i)} $d(q_1,p),d(q_2,p)<r$;
\medskip
\myitem{(ii)} $\|X(q_0)\|\leqslant 1$ and $X$ lies strictly above $TP$ at $q_0$; and
\medskip
\myitem{(iii)} $\|X(q_1)\|=1$ and $X$ lies strictly above $TP$ at $q_1$, making an angle of at least $\varphi$ with $TP$ at that point.
\medskip
\noindent Then $X$ lies strictly above $TP$ at every point of the geodesic joining $q_0$ to $q_1$.
\endproclaim
\proclabel{LemmaJacobiFieldsLieAbove}
\proof Assume the contrary. Let $(r_n)_\ninn$ be a sequence converging to $0$. For all $n$, let $p_n\in K$ be a point, $P_n$ a geodesic hyperplane at $p_n$, $q_{1,n},q_{2,n}$ two coaxial points in $P_n$ and $X_n$ a Jacobi field over the geodesic joining $q_{1,n}$ to $q_{2,n}$ such that:
\medskip
\myitem{(i)} $\opMax(d(q_{1,n},p_n),d(q_{2,n},p_n))=r_n$;
\medskip
\myitem{(ii)} $\|X_n(q_{1,n})\|\leqslant 1$ and $X_n$ lies strictly above $TP_n$ at $q_{1,n}$; and
\medskip
\myitem{(iii)} $\|X_n(q_{2,n})\|=1$ and $X_n$ lies strictly above $TP_n$ at $q_{2,n}$, making an angle of at least $\varphi$ with $TP_n$ at this point.
\medskip
\noindent Suppose, moreover, that, for all $n$, $X_n$ is tangent to $TP_n$ at some point lying between $q_{1,n}$ and $q_{2,n}$, $x_n$, say. By compactness, there exists $p_0\in K$ towards which $(p_n)_\ninn$ subconverges. Let $g$ be the Riemannian metric of $M$. For all $n$, define $g_n = r_n^{-2}g$. The sequence of pointed manifolds $(M,g_n,p_n)_\ninn$ converges towards $(\Bbb{R}^{n+1},g_{\opEuc},0)$ in the $C^\infty$ Cheeger/Gromov sense, where $g_{\opEuc}$ is the Euclidean metric over $\Bbb{R}^{n+1}$. For all $n$, $P_n$ is also a geodesic hyperplane of $(M,g_n)$ and so $(P_n,p_n)_\ninn$ subconverges in the $C^\infty$ Cheeger/Gromov sense for pointed, immersed submanifolds to a pointed, affine hyperplane $(P_0,0)$. Likewise, there exist coaxial points $q_{1,0},q_{2,0}\in P_0$, a Jacobi field $X_0$, and a point $x_0$ lying between $q_{1,0}$ and $q_{2,0}$ towards which $(q_{1,n})_\ninn$, $(q_{2,n})_\ninn$, $(r_n X_N)_\ninn$ and $x_0$ subconverge respectively. Moreover:
\medskip
\myitem{(i)} $\opMax(d(q_{1,0},0),d(d_{2,0},0)) = 1$;
\medskip
\myitem{(ii)} $\|X_0(q_{1,0})\|\leqslant 1$ and $X_0$ lies (not necessarily strictly) above $TP$ at $q_{1,0}$; and
\medskip
\myitem{(iii)} $\|X_0(q_{2,0})\|=1$ and $X_0$ lies strictly above $TP$ at $q_{2,0}$.
\medskip
\noindent It follows that $X_0$ is not tangent to $P$ at any point along the closed geodesic joining $q_{1,0}$ to $q_{2,0}$, except possibly at $q_{1,0}$. Moreover, if $X_0$ is tangent to $P$ at $q_{1,0}$, then its derivative in the direction normal to $P$ at this point is non vanishing. However, $X_0$ is tangent to $TP$ at $x_0$. It follows from the first assertion that $x_0=q_{1,0}$, but then the derivative of $X_0$ in the direction normal to $P$ at $q_{1,0}$ vanishes, and this contradicts the second assertion. This is absurd and the result follows.\qed
\medskip
\noindent This lemma allows us to prove Propostion \procref{PropClosednessOfSemiConvexity}:
\medskip
{\bf\noindent Proof of Proposition \procref{PropClosednessOfSemiConvexity}:\ }Let $K_{t_0}$ be the set bounded by $\tilde{\Sigma}_{t_0}$ and $H_{t_0}$. Let $\gamma:[0,1]\rightarrow M$ be a geodesic above $H_{t_0}$ with endpoints in $K_{t_0}$. We aim to show that the whole of $\gamma$ is contained in $K_{t_0}$. It suffices to consider the case where both endpoints of $\gamma$ lie in $\tilde{\Sigma}_{t_0}$. The remaining cases are similar and much simpler. Recall that $\tilde{\Sigma}_{t_0}$ divides into two components, $\hat{\Sigma}^-_{t_0}$ and $\Sigma_{t_0}$. These components have different properties and we thus consider the various resulting cases seperately. Let $\hat{K}_{t_0}$ be the set bounded by $\hat{\Sigma}_{t_0}$ and $H_{t_0}$. We may assume that $\hat{K}_{t_0}$ is semi-convex. Since the endpoints of $\gamma$ lie in $\tilde{\Sigma}_{t_0}\subseteq\hat{K}_{t_0}$, the whole of $\gamma$ therefore lies in $\hat{K}_{t_0}$. Thus, by choosing $t_1$ sufficiently small, we may assume that $\gamma$ is sufficiently short to satisfy the hypotheses of Proposition \procref{PropClosednessOfSemiConvexity} with $\varphi=\theta/2$.
\medskip
\noindent Suppose that $\gamma$ lies strictly above $H_{t_0}$. Then there exists $\epsilon>0$ such that $\gamma$ lies above $H_{t_0+\epsilon}$. Since $\tilde{\Sigma}_{t_0+\epsilon}$ is semi-convex, $\gamma$ lies in $K_{t_0+\epsilon}\subseteq K_{t_0}$ and the result follows in this case. We thus assume that $\gamma$ meets $H_{t_0}$ at some point, $s\in[0,1]$.
\medskip
\noindent Suppose that $\gamma$ is transverse to $H_{t_0}$ at $s$. Then, $s$ is an endpoint of $[0,1]$ and, without loss of generality, $s=0$. By strict convexity of $H_{t_0}$, $\gamma(]0,1])$ lies strictly above $H_{t_0}$. Suppose that $\gamma(0)$ lies in $\hat{\Sigma}^-_{t_0}$. By Proposition \procref{CorTransversalityOfFixedPart}, both $\hat{\Sigma}$ and $\Gamma$ are transverse to $H_{t_0}$ at this point. There thus exists a smooth curve $\eta:[0,\epsilon[\rightarrow M$ such that:
\medskip
\myitem{(i)} $\eta(0)=\gamma(0)$;
\medskip
\myitem{(ii)} $\partial_t\eta(0)$ is transverse to $TH_{t_0}$;
\medskip
\myitem{(iii)} for $s>0$, $\eta(s)$ lies strictly above $H_{t_0}$; and
\medskip
\myitem{(iv)} for all $s$, $\eta(s)$ lies in $\hat{\Sigma}^-$.
\medskip
\noindent For all $s\in[0,\epsilon[$, let $\gamma_s$ be the unique geodesic joining $\eta(s)$ to $\gamma(1)$. For sufficiently small $s$, $\gamma_s$ lies strictly above $H_{t_0}$. Since, for all $t>t_0$, $\tilde{\Sigma}_t$ is semi-convex, for all sufficiently small $s$, $\gamma_s$ is contained in $K_{t_0}$. The result follows in this case by taking limits.
\medskip
\noindent Suppose that $\gamma(0)$ lies in $\Sigma\setminus\Gamma$. By Proposition \procref{PropAngleBddBelow}, after reducing $t_1$ if necessary, the outward pointing normal to $\Sigma$ makes an angle of at least $\theta/2$ with $H_{t_0}$ at $\gamma(0)$. There therefore exists a smooth curve $\eta:[0,\epsilon[\rightarrow M$ such that:
\medskip
\myitem{(i)} $\eta(0)=\gamma(0)$;
\medskip
\myitem{(ii)} $\partial_t\eta(0)$ is transverse to $TH_{t_0}$;
\medskip
\myitem{(iii)} for $s>0$, $\eta(s)$ lies strictly above $H_{t_0}$; and
\medskip
\myitem{(iv)} for all $s$, $\eta(s)$ lies in $K_{t_0}$.
\medskip
\noindent For all $s\in[0,\epsilon[$, let $\gamma_s$ be the unique geodesic joining $\eta(s)$ to $\gamma(1)$. For sufficiently small $s$, $\gamma_s$ lies strictly above $H_{t_0}$, and the result follows in this case as before. This completes the case where $\gamma$ is transverse to $H_{t_0}$ at $s$, and we thus suppose that $\gamma$ is tangent to $H_{t_0}$ at $s$.
\medskip
\noindent Let $P$ be the geodesic hyperplane tangent to $H_{t_0}$ at $\gamma(s)$. Suppose that $\gamma(0)$ and $\gamma(1)$ both lie in $\hat{\Sigma}^-_{t_0}\setminus\Gamma$. Since $\hat{\Sigma}$ bounds a strictly convex set, $K$, $\gamma$ is transverse to $\hat{\Sigma}^-$ at $\gamma(0)$ and $\gamma(1)$ (for otherwise, by strict convexity, it could only intersect $\hat{\Sigma}^-$ at one point, which is absurd). Let $X$ be a Jacobi field over the geodesic joining $\gamma(0)$ and $\gamma(1)$ such that $X$ equals the unit upward pointing normal to $P$ at both endpoints. By Lemma \procref{LemmaJacobiFieldsLieAbove}, $X$ lies everywhere above $TP$. Thus, if $\gamma_t$ is a geodesic variation of $\gamma$ with Jacobi field $X$, then, for sufficiently small $t$, $\gamma_t$ lies strictly above $P$ and therefore also above $H_{t_0}$. Moreover, by transversality, for sufficiently small $t$, $\gamma_t$ intersects $\hat{\Sigma}^-_{t_0}$ at two points near $\gamma(0)$ and $\gamma(1)$. We thus obtain a family of geodesic segments lying strictly above $H_{t_0}$ with endpoints in $\tilde{\Sigma}_{t_0}$ converging towards $\gamma$. By semi-convexity, all these geodesic segments are contained within $K_{t_0}$, and thus, taking limits, $\gamma$ is contained within $K_{t_0}$. This proves the result in this case.
\medskip
\noindent Suppose that $\gamma(0)$ lies in $\hat{\Sigma}^-_{t_0}\setminus\Gamma$ and $\gamma(1)$ lies in $\Sigma\setminus\Gamma$. As before, $\gamma$ is transverse to $\hat{\Sigma}$ at $\gamma(0)$. By Proposition \procref{PropAngleBddBelow}, after reducing $t_1$ if necessary, the outward pointing normal to $\Sigma$ makes an angle of at least $\theta/2$ with $TP$ at $\gamma(1)$. Let $X$ be a Jacobi field over $\gamma$ such that $X(0)$ is the upward pointing normal vector over $P$ at $\gamma(0)$ and $X(1)$ points into $K_{t_0}$ making an angle of at least $\theta/2$ with $TP$ at $\gamma(1)$. By Lemma \procref{LemmaJacobiFieldsLieAbove}, $X$ lies everywhere above $TP$. Thus, if $\gamma_t$ is a geodesic variation of $\gamma$ with Jacobi field $X$, then, for sufficiently small $t$, $\gamma_t$ lies strictly above $P$ and therefore also above $H_{t_0}$. Moreover, for small $t$, $\gamma_t(1)$ lies inside $K_{t_0}$, and, by transversality, $\gamma_t$ intersects $\hat{\Sigma}^-_{t_0}$ at some point near $\gamma(0)$. We thus obtain a family of geodesic segments lying strictly above $H_{t_0}$ with endpoints in $K_{t_0}$ converging towards $\gamma$. By semi-convexity, all these geodesic segments are contained within $K_{t_0}$, and thus, taking limits, $\gamma$ is contained within $K_{t_0}$. This proves the result in this case.
\medskip
\noindent Suppose that both $\gamma(0)$ and $\gamma(1)$ lie in $\Sigma\setminus\Gamma$. By Proposition \procref{PropAngleBddBelow}, after reducing $t_1$ if necessary, the outward pointing normal to $\Sigma$ makes an angle of at least $\theta/2$ with $P$ at both these points. Let $X$ be a Jacobi field over $\gamma$ such that both $X(0)$ and $X(1)$ point into $K_{t_0}$ at $\gamma(0)$ and $\gamma(1)$ respectively, making an angle of at least $\theta/2$ with $TP$ at these points. By Lemma \procref{LemmaJacobiFieldsLieAbove}, $X$ lies everywhere above $TP$, and the result follows in this case as before.
\medskip
\noindent We now consider the case where at least one end point of $\gamma$ lies on $\Gamma$. Suppose that $\gamma(0)$ lies on $\Gamma$ but $\gamma(1)$ doesn't. By Proposition \procref{CorTransversalityOfFixedPart}, $\Gamma$ is transverse to $P$ at $\gamma(0)$. Let $X$ be a Jacobi field over $\gamma$ such that $X(0)$ is tangent to $\Gamma$ and points strictly upwards from $P$ at $\gamma(0)$. If $\gamma(1)$ lies in $\hat{\Sigma}^-_{t_0}$, then we suppose that $X(1)$ is the upward pointing unit normal over $P$ at $\gamma(1)$. If $\gamma(1)$ lies in $\Sigma_{t_0}$, then we assume that $X(1)$ points into $K_{t_0}$ at $\gamma(1)$, making an angle of at least $\theta/2$ with $TP$ at this point. By Lemma \procref{LemmaJacobiFieldsLieAbove}, $X$ lies everywhere above $TP$, and the result follows in this case as before.
\medskip
\noindent Finally suppose that both $\gamma(0)$ and $\gamma(1)$ lie on $\Gamma$. It follows by Proposition \procref{PropGeodesicsWithBothEndsInBoundary} that, after increasing $t_1$ if necessary, there exists a small deformation of $\gamma$ whose end points both lie on $\hat{\Sigma}_{t_0}^-\setminus\Gamma$. We thus reduce this case to an earlier case, and this completes the proof.\qed
\goodbreak
\newhead{Immersed Boundaries}
\noindent Let $M^{n+1}$ be an $(n+1)$-dimensional manifold. We recall that the reasoning of Section \headref{FirstOrderUpperBounds} is only valid when the boundary is embedded. We now show how this reasoning may be adapted by a simple modification to also treat the case where the boundary is permitted to have self intersections.
\medskip
\noindent Let $\Gamma^{n-1}=(i,(G^{n-1},\partial G^{n-1}))$ be a compact, codimension $2$, immersed submanifold in $M$. We say that $\Gamma$ is {\bf generic} if and only if, for all $p\neq q$ such that $i(p)=i(q)$:
\headlabel{HeadImmersedBoundaries}
$$
T_p\Gamma \neq T_q\Gamma.
$$
\noindent This definition is motivated by the following elementary result:
\proclaim{Proposition \nextprocno}
\myitem{(i)} Let $\Gamma\subseteq M$ be a compact, codimension $2$, immersed submanifold. There exists a sequence $(\Gamma_n)_\ninn$ of generic, compact, codimension $2$, immersed submanifolds which converges to $\Gamma$ in the $C^\infty$ sense.
\medskip
\myitem{(ii)} Let $(\Gamma_t)_{t\in[0,1]}\subseteq M$ be a smooth family of compact, codimension $2$, immersed submanifolds such that $\Gamma_0$ and $\Gamma_1$ are generic. There exists a sequence $(\Gamma_{n,t})_{\ninn}$ of smooth families of generic, compact, codimension $2$, immersed submanifolds such that:
\medskip
\myitem{(a)} for all $n$, $\Gamma_{n,0}=\Gamma_0$ and $\Gamma_{n,1}=\Gamma_1$; and
\medskip
\myitem{(b)} $(\Gamma_{n,t})_{\ninn}$ converges to $(\Gamma_t)$ in the $C^\infty$ sense.
\endproclaim
\proclabel{PropDensityOfGenericity}
\proof This follows from Sard's Lemma in the usual manner. Explicitely, a generic codimension $2$ immersion self-intersects over a submanifold of codimension $4$, from which $(i)$ follows, and every immersion in a generic isotopy of codimension $2$ immersions self-intersects over a submanifold of codimension $3$, from which $(ii)$ follows. See \cite{GuillPoll} for details.\qed
\medskip
\noindent Let $(\Gamma_n)_\ninn$ be a sequence of strictly convex, codimension $2$, immersed submanifolds with convexity orientation. For all $n\in\Bbb{N}\munion\left\{0\right\}$, let $\msf{N}^+_n$ be the convexity coorientation of $\Gamma_n$. Suppose that $(\Gamma_n)_\ninn$ converges in the $C^\infty$ sense to a strictly convex, codimension $2$, immersed submanifold, $\Gamma_0$ and suppose, moreover, that $\Gamma_0$ is generic. In particular, by taking a subsequence, we may suppose that $\Gamma_n$ is also generic for all $n$.
\proclaim{Lemma \nextprocno}
\noindent Choose $\theta>0$. There exists $r>0$ such that if $(\Sigma_n)_\ninn$ is a sequence of strictly convex, immersed hypersurfaces such that, for all $n$:
\medskip
\myitem{(i)} $\partial\Sigma_n=\Gamma_n$; and
\medskip
\myitem{(ii)} the outward pointing unit normal over $\Sigma_n$ makes an angle of at least $\theta$ with $\msf{N}^+_n$ along $\Gamma_n$,
\medskip
\noindent then, for all $n$, and for all $p\in\Gamma_n$:
\medskip
\myitem{(i)} the connected component of $\Sigma_n\minter B_r(p)$ is embedded and lies on the boundary of a convex subset of $B_r(p)$; and
\medskip
\myitem{(ii)} this connected component only meets one connected component of $\Gamma_n\minter B_r(p)$.
\endproclaim
\proclabel{LemmaBoundaryCompactnessWhenNotEmbedded}
\remark Using this result in conjunction with the compactness of the family of bounded convex sets, we obtain
$C^{0,\alpha}$ compactness near the boundary for families of locally convex immersed hypersurfaces. In particular, this result may be used to extend the conclusions of \cite{SmiNLD} to the case of compact hypersurfaces with non-trivial boundary (see \cite{SmiGPP}).
\medskip
\proof For all $n\in\Bbb{N}\munion\left\{0\right\}$, choose $p_n\in\Gamma_n$ and suppose that $(p_n)_\ninn$ converges to $p_0$. For all $n\in\Bbb{N}\munion\left\{0\right\}$, let $q_n\in M$ be the image of $p_n$. Choose $r>0$ such that, for all $n\in\Bbb{N}\munion\left\{0\right\}$, the connected component of $\Gamma_n\minter B_r(q_n)$ containing $p_n$ is embedded, and denote this component by $\Gamma_{n,0}$. For all $n$, we identify $M$ with $B_r(p_n)$, reducing $r$ whenever necessary.
\medskip
\noindent As in Section \headref{FirstOrderUpperBounds}, for all $n\in\Bbb{N}\munion\left\{0\right\}$, let $H_n$ be a strictly convex, embedded hypersurface tangent to $\Gamma_n$ at $p_n$ such that:
\medskip
\myitem{(i)} the outward pointing normal to $H_n$ at $p_n$ makes an angle of no more than $\theta/2$ with $\msf{N}^+_n$ at $p_n$; and
\medskip
\myitem{(ii)} the shape operator of $H_n$ is everywhere strictly bounded above by $\delta\opId$, where $\delta$ is small.
\medskip
\noindent We suppose, moreover, that $(H_n)_\ninn$ converges to $H_0$ in the $C^\infty$ sense. Likewise, as in Section \headref{FirstOrderUpperBounds}, for all $n\in\Bbb{N}$, we extend $H_n$ to a foliation $(H_{n,t})_{t\in\Bbb{R}}$.
\medskip
\noindent Since $\Gamma_0$ is generic, we may suppose that $H_0$ is transverse at $q_0$ to every connected component of $\Gamma_0\minter B_r(q_0)$ not equal to $\Gamma_{0,0}$ which passes through $q_0$. Thus, reducing $r$ if necessary, for all $n$, if $\Gamma'_{n,0}$ is a connected component of $\Gamma_n\minter B_r(q_n)$ which is different from $\Gamma_{n,0}$, then $\Gamma'_{n,0}$ is transverse to $H_{n,t}$, for all $t$.
\medskip
\noindent Let $t_0<0$ be as in Section \headref{FirstOrderUpperBounds}, and, for all $t\in]t_0,0[$, let $\Sigma_{n,t}$ be the connected component of $\Sigma_n$ containing $p_n$ which lies above $H_{n,t}$. Define $T$ to be the set of all $t\in]-t_0,0[$ such that $\Gamma_{n,0}$ is the only connected component of $\Gamma_n\minter B_r(q_n)$ which intersects $\Sigma_{n,t}$. Trivially, $T$ is non-empty. Let $t_1=\minf T$ and suppose that $t_1>t_0$. Let $\Gamma'_{n,0}\neq\Gamma_{n,0}$ be the connected component of $\Gamma_n\minter B_r(q_n)$ which intersects $\Sigma_{n,t_1}$. For $t>t_1$, the reasoning of Section \headref{FirstOrderUpperBounds} proceeds as in the case where the boundary is embedded, and it follows that $\Sigma_{n,t_1}$ is embedded, is transverse to $H_{t_1}$ and bounds a semi-convex set above $H_{t_1}$. $\Gamma'_{n,0}$ is therefore tangent to $H_{n,t_1}$ at the point of intersection, since, otherwise $\Gamma'_{n,0}$ would intersect $\Sigma_{n,t}$ non trivially at some point lying above $H_{n,t_1}$, which is absurd. However, this contradicts the definition of $r$. It follows that $t_1=t_0$, and the result now follows as in the case of Lemma \procref{LemmaBoundaryCompactnessOfConvexImmersions} by taking intersections with a ball of radius less than $t_0$.\qed
\goodbreak
\newhead{First Order Lower Bounds}
\noindent Let $M^{n+1}$ be an $(n+1)$-dimensional Riemannian manifold. Let $\Gamma^{n-1}\subseteq M$ be a generic, strictly convex, codimension $2$, immersed submanifold with convexity orientation. Let $A_\Gamma$ be the shape operator of $\Gamma$ and let $\msf{N}^-$ and $\msf{N}^+$ be the convexity orientation and coorientation respectively of $\Gamma$. As in \cite{CaffNirSprI}, second order bounds require uniform lower bounds on the angle between $\msf{N}^-$ and the normal to any hypersurface of constant Gaussian curvature with boundary equal to $\Gamma$. This is guaranteed by the following result:
\proclaim{Proposition \nextprocno}
\noindent For all $k>0$, there exists $\phi>0$ (which only depends on $M$, $\Gamma$ and $\theta$) such that if $(\Sigma^n,\partial\Sigma^n)$ is a smooth, convex immersed hypersurface such that:
\medskip
\myitem{(i)} $\partial\Sigma=\Gamma$;
\medskip
\myitem{(ii)} the Gaussian curvature of $\Sigma$ is at least $k$; and
\medskip
\myitem{(iii)} the outward pointing normal to $\Sigma$ over $\Gamma$ makes an angle of at least $\theta$ with $\msf{N}^+(p)$,
\medskip
\noindent then the outward pointing normal to $\Sigma$ over $\Gamma$ also makes an angle of at least $\phi$ with $\msf{N}^-(p)$.
\endproclaim
\proclabel{PropFirstOrderLowerBound}
\noindent Let $r>0$ and let $\Sigma$ be a $C^{0,1}$ locally convex hypersurface in $M$ such that:
\medskip
\myitem{(i)} $\partial\Sigma\subseteq\Gamma\munion B_r(p)$;
\medskip
\myitem{(ii)} $\Sigma$ is compatible with the orientation on $\Gamma$;
\medskip
\myitem{(iii)} the outward pointing normal to $\Sigma$ along $\Gamma$ always makes an angle of at least $\theta$ with $\msf{N}^+$; and
\medskip
\myitem{(iv)} the outward pointing normal to $\Sigma$ at $p$ coincides with $\msf{N}^-(p)$.
\medskip
\noindent Let $\opSymm(\Bbb{R}^n)$ denote the set of positive definite, symmetric matrices over $\Bbb{R}^n$. For $t>0$, we define $F_t\subseteq\opSymm(\Bbb{R}^n)$ by:
$$
F_t = \left\{ A\in\opSymm(\Bbb{R}^n)\text{ s.t. }A\geqslant 0\ \&\ \opDet(A)\geqslant t\right\}.
$$
\noindent Observe that if $A\in F_t$ and if $M\geqslant 0$, then $A+M\in F_t$. In the language of \cite{HarveyLawsonI}, this implies that $F_t$ is a Dirichlet set. In particular, if $A\notin F_t$ and $M\geqslant 0$, then $A-M\notin F_t$. Proposition \procref{PropFirstOrderLowerBound} is proven using barriers, which are constructed using the following result:
\proclaim{Proposition \nextprocno}
\noindent Choose $\delta>0$. There exists a neighbourhood $U$ of $p$ and a smooth function $f:U\rightarrow\Bbb{R}$ such that:
\medskip
\myitem{(i)} $f\geqslant 0$ along $\partial (U\minter\Sigma)$;
\medskip
\myitem{(ii)} there exists $q\in U\minter\Sigma$ such that $f(q)<0$; and
\medskip
\myitem{(iii)} for all $q\in B_r(p)$, the shape operator of the level subset of $f$ passing through $q$ with respect to $\nabla f$ is conjugate to an element of $F_\delta^c$.
\endproclaim
\proclabel{PropConstructionOfBarrier}
\noindent Let $S$ be a smooth, immersed hypersurface in $M$ such that:
\medskip
\myitem{(i)} $\partial S=\partial\Gamma$;
\medskip
\myitem{(ii)} the upward pointing normal to $S$ at $p$ is equal to $\msf{N}^-(p)$; and
\medskip
\myitem{(iii)} the shape operator of $S$ at $p$ is supported along the subspace $T_p\Gamma$.
\medskip
\noindent Let $H$ be a strictly concave immersed hypersurface in $M$ such that:
\medskip
\myitem{(i)} the downward pointing normal to $H$ at $p$ lies in $X_p$ and makes an angle of at most $\theta/2$ with $\msf{N}^+(p)$; and
\medskip
\myitem{(ii)} $\Gamma$, $\Sigma$ and $S$ locally lie strictly above $H$.
\medskip
\noindent Let $d_p$, $d_S$ and $d_H$ denote the (signed) distance in $M$ to $p$, $S$ and $H$ respectively. Observe that $(\nabla d_S,\nabla d_H)$ is a linearly independant pair which spans the space of normal vectors to $\Gamma$ at $p$. For any two functions, $f$ and $g$, we define the $(n-2)$-dimensional distribution, $E(f,g)$, near $p$ by:
$$
E(f,g) = \langle \nabla f,\nabla g\rangle^\perp,
$$
\noindent where $\langle U,V\rangle$ here represents the subspace spanned by the vectors $U$ and $V$. Let $e_1,...,e_{n-1}$ be an orthonormal basis for $T_p\Gamma$ with respect to which $A_\Gamma(\msf{N}^-)$ is diagonal. Let $\lambda_1,...,\lambda_{n-1}$ be the corresponding eigenvalues. We may suppose that $0=\lambda_1\leqslant\lambda_2\leqslant ...\leqslant \lambda_{n-1}$. We extend $(e_1,...,e_{n-1})$ to a local frame in $TM$ such that, for all vectors, $X$, at $p$:
$$\matrix
\langle\nabla_X e_i,\nabla d_S\rangle = -\opHess(d_S)(X,e_i),\hfill\cr
\langle\nabla_X e_i,\nabla d_H\rangle = -\opHess(d_H)(X,e_i).\hfill\cr
\endmatrix$$
\noindent Define the distribution $E$ near $P$ to be the span of $e_1,...,e_{n-1}$.
\proclaim{Proposition \nextprocno}
\noindent If $D$ represents the Grassmannian distance between two $(n-1)$-dimensional subspaces, then:
$$
D(E,E(d_S,d_H)) = O(d_p^2).
$$
\endproclaim
\proclabel{PropDistanceInGrassmannianA}
\proof By definition of $e_i$, for all vectors $X$ at $p$:
$$
X\langle e_i,\nabla d_S\rangle = X\langle e_i,\nabla d_H\rangle = 0.
$$
\noindent The result follows.\qed
\medskip
\noindent For any smooth function, $f$, we define $D(f,E)$ by:
$$
D(f,E) = \opDet(\opHess(f)|_E),
$$
\noindent where $\opHess(f)|_E$ is the restriction of the Hessian of $f$ to $E$.
\proclaim{Proposition \nextprocno}
\noindent Let $f$ be such that $f(p),\nabla f(p)=0$ and the restriction of $\opHess(f)$ to $H$ at $p$ is positive definite. There exists a function $x$ such that $x(p),\opHess(x)(p)=0$ and:
$$
D(d_S + x(d_H - f),E) = O(d_p)^2.
$$
\endproclaim
\proclabel{PropCorrectionTerm}
\proof The Hessian of $xf$ vanishes at $p$. Likewise, the Hessian of the second order term $xd_H$ vanishes over $(\nabla d_H)^\perp$ and therefore over $E$ at $p$. It follows that the term $x(d_H -f)$ does not affect the restriction of the Hessian of the function to $E$ at $p$. Thus:
$$
\nabla D = \opTr(\opAdj(\opHess(d_S)|_E)\nabla(\opHess(d_S + x(d_H - f))(e_i,e_j))),
$$
\noindent where $\opAdj(\opHess(d_S)|_E)$ is the adjugate matrix of $\opHess(d_S)|_E$. If more than one of the eigenvalues of $\opHess(d_S)|_E$ vanishes, then $\opAdj(\opHess(d_S)|_E)$ also vanishes, and the result follows trivially by taking $x=0$. Suppose therefore that only one eigenvalue of $\opHess(d_S)|_E$ vanishes. Let $\mu_1,...,\mu_{n-1}$ be the eigenvalues of the adjugate matrix, then $\mu_1=\lambda_2...\lambda_{n-1}$ and $\mu_2=...=\mu_{n-1}=0$. Define the vectors $U$ and $V$ at $p$ by:
$$\matrix
U \hfill&= \nabla D(d_S,E),\hfill\cr
V \hfill&= \nabla D(d_S + x(d_H-f),E).\hfill\cr
\endmatrix$$
\noindent Denote $P=x(d_H-f)$. At $p$:
$$
\opHess(P) = \nabla x\otimes \nabla d_H + \nabla d_H\otimes\nabla x.
$$
\noindent At $p$, for all $i$, by definition, $\langle e_i,\nabla d_H\rangle = 0$. Thus, recalling the formula for $\nabla e_i$:
$$\matrix
X\opHess(P)(e_i,e_j) \hfill&= (\nabla_X\opHess(P))(e_i,e_j) + \opHess(P)(\nabla_X e_i,e_j) + \opHess(P)(e_i,\nabla_X e_j)\hfill\cr
&=(\nabla_X\opHess(P))(e_i,e_j)\hfill\cr
&\qquad\qquad + \langle\nabla x, e_j\rangle\langle\nabla_X e_i,\nabla d_H\rangle
+ \langle\nabla x, e_i\rangle\langle\nabla_X e_j,\nabla d_H\rangle\hfill\cr
&=(\nabla_X\opHess(P))(e_i,e_j) - \opHess(d_H)(X,e_i)x_{;j} - \opHess(d_H)(X,e_j)x_{;i}.\hfill\cr
\endmatrix$$
\noindent We extend $(e_i)_{1\leqslant i\leqslant n-1}$ to an orthonormal basis $(e_i)_{0\leqslant i\leqslant n}$ for $T_pM$. With respect to this basis, for all $k$:
$$
\frac{1}{\mu_1}\langle V - U,e_k\rangle = (d_{H;11}-f_{;11})x_{;k} - 2f_{;1k}x_{;1}.
$$
\noindent Consider the linear map, $M$, given by:
$$
(M\xi)_{k} = (d_{H;11}-f_{;11})\xi_k - 2f_{;1k}\xi_1.
$$
\noindent Suppose that $M\xi=0$. Then, in particular, bearing in mind that $d_{H;11}\leqslant 0$ and $f_{;11}\geqslant 0$:
$$\matrix
&(d_{H;11}-3f_{;11})\xi_{1} \hfill&= 0\hfill\cr
\Rightarrow\hfill& \xi_1 \hfill&= 0\hfill\cr
\Rightarrow\hfill& \xi\hfill&=0.\hfill\cr
\endmatrix$$
\noindent $M$ is therefore invertible, and there exists $\xi$ such that:
$$
M\xi = -U.
$$
\noindent If we define $x$ such that:
$$
x(p) = 0,\qquad \nabla x(p) = \xi, \qquad \opHess(x)(p) = 0,
$$
\noindent then:
$$
\nabla D(d_S + x(d_H - f),E) = 0.
$$
\noindent This completes the proof.\qed
\medskip
\noindent Define $\Phi_0$ by:
$$
\Phi_0 = d_S + x(d_H - f).
$$
\noindent For $M>0$, define $\Phi$ by:
$$
\Phi = d_S + x(d_H - f) + Md_H^2.
$$
\proclaim{Proposition \nextprocno}
\noindent If $D$ represents the Grassmannian distance between two $(n-2)$-dimensional subspaces then:
$$
D(E(d_S,d_H),E(\Phi,d_H)) = O(d_p^2) + O(d_H).
$$
\endproclaim
\proof Since $xf$ is of order $3$ at $p$:
$$
\nabla\Phi = \nabla d_S + (x+2Md_H)\nabla d_H + O(d_p^2) + O(d_H).
$$
\noindent Thus:
$$
\langle\nabla\Phi,\nabla d_H\rangle = \langle\nabla d_S + O(d_p^2) + O(d_H),\nabla d_H\rangle,
$$
\noindent where $\langle\cdot,\cdot\rangle$ here represents the subspace generated by two vectors. The result follows.\qed
\proclaim{Corollary \nextprocno}
\noindent If $D$ represents the Grassmannian distance between two $(n-2)$-dimensional subspaces, then:
$$
D(E,E(\Phi,d_H)) = O(d_p^2) + O(d_H).
$$
\endproclaim
\proclabel{CorGrassmannianDistanceB}
\proof This follows from the triangle inequality and Proposition \procref{PropDistanceInGrassmannianA}.\qed
\medskip
\noindent Finally, we recall the following technical property of convex sets. Let $UM$ be the bundle of unit spheres in $TM$. Let $K\subseteq M$ be a compact, convex set with non-trivial interior. For all $q\in\partial K$, let $\Cal{N}(q)\subseteq U_qM$ be the set of supporting normals to $K$ at $q$. This set is a closed, convex subset subset of $U_qM$. Moreover, we have the following continuity result:
\proclaim{Proposition \nextprocno}
\noindent Let $q_0,(q_n)_\ninn\in\partial K$ be such that $(q_n)_\ninn$ converges to $q_0$. For all $n$, let $N_n$ be an element of $\Cal{N}(q_n)$. If $d$ denotes the distance in $UM$, then $(d(N_n,\Cal{N}(q_0)))_\ninn$ converges to $0$.
\endproclaim
\proclabel{PropContinuityOfNormal}
\noindent We now prove Proposition \procref{PropConstructionOfBarrier}:
\medskip
{\noindent\bf Proof of Proposition \procref{PropConstructionOfBarrier}:\ } For $\epsilon>0$, define the open set $U_\epsilon\subseteq M$ by:
$$
U_\epsilon = \left\{p\in M\text{ s.t. }d_p(x)<\epsilon\text{ and }d_H(x)<\epsilon^2\right\}.
$$
\noindent $\partial(\Sigma\minter U_\epsilon)$ consists of two components: $\partial\Sigma\minter U_\epsilon = \Gamma\minter U_\epsilon$ and $\partial U_\epsilon\minter\Sigma$. We first obtain lower estimates for $\Phi_0$ along these two components.
\medskip
\noindent We choose $f$ such that, along $\Gamma$, $(f-d_H)=O(d_p^3)$. Consequently, $x(f-d_H)=O(d_p^4)$ along $\Gamma$. Thus, since $O(d_p^2)=O(d_H)$ along $\Gamma$ and since $d_S$ vanishes along $\Gamma$, there exists $K_1>0$ such that, along $\Gamma$:
$$
\left|d_S + x(d_H - f)\right|\leqslant K_1d_H^2.
$$
\noindent This yields lower bounds for $\Phi_0$ along $\partial\Sigma\minter U_\epsilon$.
\medskip
\noindent Since $\Sigma$ is a convex immersion, and since $\partial\Sigma=\Gamma$ is smooth, $\Sigma$ has a unique supporting normal at $p$, which coincides with $\nabla d_S$. Now let $V$ be a field of unit vectors defined near $p$ such that $V(p)$ makes an angle of exactly $\theta$ with $\msf{N}^+(p)$. For $q\in M$, let $U_qM$ be the unit sphere in $T_qM$. Let $D_q$ be the distance in $U_qM$ and let $C_q$ be the shortest geodesic in $U_qM$ joining $V(q)$ to $\nabla d_S(q)$. Near $p$, $V(p)$, $\nabla d_S$ and $-\nabla d_{H}$ are configured as shown in Figure \figref{FigConfiguration}:
\placefigure{}{%
\placelabel[1.9][1.85]{$\nabla d_S$}%
\placelabel[5.65][0.95]{$-\nabla d_{H}$}%
\placelabel[1.9][3.2]{$C$}%
\placelabel[5.55][1.85]{$V$}%
}{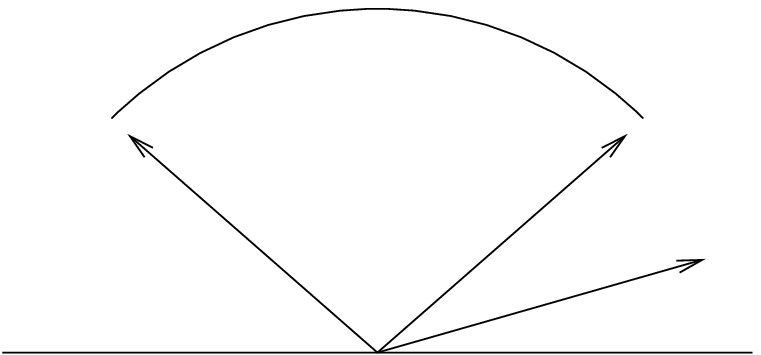}{}
\figlabel{FigConfiguration}
\noindent However, by definition, $\nabla d_S(p)=\msf{N}_p$. By Property $(iii)$ of $\Sigma$, we can extend $\Sigma$ slightly beyond $\Gamma$ to a $C^{0,1}$ locally strictly convex hypersurface whose set of supporting normals at $p$ is contained within $C_p$. Thus, by Proposition \procref{PropContinuityOfNormal}, there exists a continuous function $\delta:[0,\infty[\rightarrow[0,\infty[$ such that $\delta(0)=0$ and, for all $q\in\Sigma$, if $\msf{N}_q$ is a supporting normal to $\Sigma$ at $q$, then:
$$
D_q(\msf{N}_q,C_q)\leqslant\delta(d_p(q)).
$$
\noindent Thus, if, for all $q\in\Sigma$, $\pi_q$ is the orthogonal projection onto a supporting hyperplane of $\Sigma$ at $q$, then:
\medskip
\myitem{(i)} there exists $c>0$ such that, for all $q$ sufficiently close to $p$:
$$
\|\pi_1(\nabla d_{H})\|\geqslant c;\text{ and}
$$
\myitem{(ii)} for all $q$ sufficiently close to $p$:
$$
\langle\pi_q(\nabla d_S),\pi_q(\nabla d_{H})\rangle \geqslant -\delta(d_p(q)).
$$
\noindent Now consider $q_0\in\Sigma\minter\partial U$. Let $\gamma:I\rightarrow\Sigma$ be an integral curve of $\pi_q(\nabla d_{H})$ such that $\gamma(0)\in\partial\Sigma$ and $\gamma(1)=q_0$ (which is defined by approximating $\Sigma$ by smooth hypersurfaces). Bearing in mind that $d_H\geqslant 0$ along $\partial\Sigma$ and $d_S$ vanishes along $\partial\Sigma$:
$$\matrix
&d_{H}(q_0)\hfill&\leqslant \epsilon^2\hfill\cr
\Rightarrow\hfill&\opLength(\gamma)\hfill&\leqslant \epsilon^2c^{-1}\hfill\cr
\Rightarrow\hfill&(d_S\circ\gamma)(1)\hfill&\geqslant -\delta(\epsilon)\epsilon^2c^{-1}.\hfill\cr
\endmatrix$$
\noindent Thus:
$$
[d_S + x(d_{H} - f)](q_0) \geqslant -\delta(\epsilon)O(\epsilon^2),
$$
\noindent for all appropriate functions $f$ and $x$. There thus exists $\delta_1>0$ such that, along $\Sigma\minter\partial U_\epsilon$:
$$
\Phi_0> -\delta_1 d_H.
$$
\noindent Moreover, $\delta_1$ tends to $0$ as $\epsilon$ tends to $0$.
\medskip
\noindent Thus, if we choose $M=\opMax(\delta_1\epsilon^{-2},K_1)$, then $\Phi\geqslant 0$ along $\partial(\Sigma\minter U)$. Since $\opHess(\Phi_0)$ is bounded, by Proposition \procref{PropCorrectionTerm} and Corollary \procref{CorGrassmannianDistanceB}:
$$
D(\Phi_0,E(\Phi,d_H)) = O(\epsilon^2).
$$
\noindent However:
$$
\opHess(\Phi) = \opHess(\Phi_0) + 2M\nabla d_H\otimes\nabla d_H + 2Md_H\opHess(d_H).
$$
\noindent Denote:
$$
A = \frac{1}{\|\nabla\Phi\|}\opHess(\Phi)|_{\nabla\Phi^\perp}.
$$
\noindent $A$ is the shape operator of the level sets of $\Phi$. If $A$ is not non-negative definite, then it trivially lies in $F_\delta^c$. Suppose, therefore, that $A$ is non-negative definite. Let $0\leqslant\lambda_1\leqslant...\leqslant\lambda_n$ be the eigenvalues of $A$, and let $0\leqslant\lambda'_1\leqslant...\leqslant\lambda'_{n-1}$ be the eigenvalues of the restriction of $A$ to $E(\Phi,d_H)$. Observe that $\nabla d_H\otimes \nabla d_H$ vanishes on $E(\Phi,d_H)$. Moreover, since $H$ is concave, $2Md_H\opHess(d_H)$ is negative definite. It follows that the eigenvalues of the restriction of $A$ to $E(\Phi,d_H)$ are less than the eigenvalues of the restriction of $\|\nabla\Phi\|^{-1}\opHess(\Phi_0)$ to this subspace. Thus, since $\|\nabla\Phi\|$ also remains uniformly bounded away from $0$, by the preceeding calculations:
$$
\lambda'_1\cdot...\cdot\lambda'_{n-1} = O(\epsilon^2).
$$
\noindent However, by the minimax principal, for $1\leqslant i\leqslant (n-1)$:
$$
0 \leqslant \lambda_i \leqslant \lambda'_i.
$$
\noindent Thus:
$$
\lambda_1\cdot...\cdot\lambda_{n-1} = O(\epsilon^2).
$$
\noindent Consequently, since $\lambda_n = O(M) = O(\delta_1\epsilon^{-2})$, for $\epsilon$ sufficiently small:
$$
\opDet(A) < \delta.
$$
\noindent Thus $A\in F_\delta^c$, and property $(iii)$ now follows. Since $f$ is non-negative over $\partial(\Sigma\minter U_\epsilon)$, property $(i)$ also follows. Since $f(p)=0$ and $(\nabla f)(p)=\msf{N}^+(p)$, deforming $f$ slightly yields a function which still satisfies conditions $(i)$ and $(iii)$ but also satisfies condition $(ii)$. This completes the proof.\qed
\medskip
\noindent We now obtain Proposition \procref{PropFirstOrderLowerBound}:
\medskip
{\bf\noindent Proof of Proposition \procref{PropFirstOrderLowerBound}:\ }Assume the contrary. Let $(\Sigma_n,\partial\Sigma_n)_\ninn$ be a sequence of convex immersed hypersurfaces such that:
\medskip
\myitem{(i)} $\partial\Sigma_n=\Gamma$; and
\medskip
\myitem{(ii)} the Gaussian curvature of $\Sigma$ is at least $k$.
\medskip
\noindent Suppose, moreover, that there exists $(p_n)_\ninn,p_0\in\Gamma$ such that $(p_n)_\ninn$ converges to $p_0$ and the angle that the exterior normal of $\Sigma_n$ makes with $\msf{N}^-(p_n)$ at $p_n$ tends to $0$.
\medskip
\noindent By Lemma \procref{LemmaBoundaryCompactnessWhenNotEmbedded}, there exists $r>0$ such that, for all $n$, the connected component of $\Sigma_n\minter B_r(p_n)$ containing $p_n$ is embedded and bounds a convex set. For all $n$, we denote this connected component by $\Sigma_{n,0}$. By compactness of the family of convex sets, there exists a convex immersion $\Sigma_0$ to which $(\Sigma_{n,0})_\ninn$ converges in the $C^{0,\alpha}$ sense for all $\alpha$. Let $f$ be as in Proposition \procref{PropConstructionOfBarrier} with $\delta<k$. For sufficiently large $n$, $f$ achieves a strict local minimum at some interior point $q_n\in\Sigma_{n,0}$.
\medskip
\noindent Let $\opHess^0(f)$ be the Hessian of $f$ over $M$, and, for all $n$, let $\opHess^n(f)$ be the Hessian of the restriction of $f$ to $\Sigma_n$. At $q_n$:
$$
\opHess^n(f) = \opHess^0(f)|_{\nabla f^\perp} - \|\nabla f\|A_{n,0},
$$
\noindent where $A_{n,0}$ is the shape operator of $\Sigma_{n,0}$ at $q_n$. By the Maximum Principal, at $q_n$:
$$\matrix
&\opHess^0(f)|_{\nabla f^\perp} - \|\nabla f\|A_{n,0} \hfill&\geqslant 0\hfill\cr
\Rightarrow\hfill&\opHess^0(f)|_{\nabla f^\perp}\hfill&\geqslant \|\nabla f\|A_{n,0}\hfill\cr
\Rightarrow\hfill&\frac{1}{\|\nabla f\|}\opHess^0(f)|_{\nabla f^\perp}\hfill&\in F_k.\hfill\cr
\endmatrix$$
\noindent This is absurd by definition of $f$, and the result follows.\qed
\goodbreak
\newhead{Compactness}
\noindent Let $M^{n+1}$ be a Hadamard manifold. Let $(\Gamma^{n-1}_m)_\minn,\Gamma^{n-1}_0\subseteq M$ be generic, locally strictly convex, codimension $2$, immersed submanifolds with convexity orientation such that $(\Gamma_m)_\minn$ converges to $\Gamma_0$. For all $m$, let $\msf{N}^-_m$ and $\msf{N}^+_m$ be the convexity orientation and coorientation respectively of $\Gamma_m$. Let $(\phi_m)_\minn,\phi_0:M\rightarrow]0,\infty[$ be smooth, positive functions such that $(\phi_m)_\minn$ converges to $\phi_0$ in the $C_\oploc^\infty$ sense. Let $(\Sigma^n_m)_\minn\subseteq M$ be smooth, immersed, strictly convex, compact hypersurfaces such that, for all $m$:
\medskip
\myitem{(i)} $\partial\Sigma_m=\Gamma_m$;
\medskip
\myitem{(ii)} $\Sigma_m$ is compatible with the orientation of $\Gamma_m$; and
\medskip
\myitem{(iii)} the Gaussian curvature of $\Sigma_m$ at any point $p\in\Sigma_m$ is equal to $\phi_m(p)$.
\medskip
\noindent We obtain the folllowing precompactness result:
\goodbreak
\proclaim{Lemma \nextprocno}
\noindent Let $\theta\in]0,\pi[$ be an angle and let $D>0$ be a positive real number. Suppose that, for all $m$:
\medskip
\myitem{(i)} the outward pointing normal to $\Sigma_m$ makes an angle of at least $\theta$ with $\msf{N}^+_m$ at every point of $\Gamma_m$; and
\medskip
\myitem{(ii)} the diameter of $\Sigma_m$ is no greater than $D$.
\medskip
\noindent Then there exists a strictly convex, smooth immersed hypersurface, $(\Sigma_0,\partial\Sigma_0)\subseteq M$ towards which $(\Sigma_m)_\minn$ subconverges. Moreover:
\medskip
\myitem{(i)} $\partial\Sigma_0=\Gamma_0$; and
\medskip
\myitem{(ii)} the Gaussian curvature of $\Sigma_0$ at any point $p\in\Sigma_0$ is equal to $\phi_0(p)$.
\endproclaim
\proclabel{LemmaPreCompactness}
\proof By the Arzela-Ascoli Theorem of \cite{SmiAAT}, it suffices to obtain a-priori bounds for all the derivatives of the shape operators of the hypersurfaces $(\Sigma_m)_\minn$. For all $m$, let $A_m$ be the shape operator of $\Sigma_m$. Let $(p_m)_\minn,p_0$ be points such that:
\medskip
\myitem{(i)} for all $m$, $p_m\in\Gamma_m$; and
\medskip
\myitem{(ii)} $(p_m)_\minn$ converges to $p_0$.
\medskip
\noindent Choose $\epsilon>0$. There exists $r_1>0$ and, for all $m$, a smooth, embedded, strictly locally convex hypersurface $\hat{\Sigma}_m$ such that:
\medskip
\myitem{(i)} $p_m\in\hat{\Sigma}_m$;
\medskip
\myitem{(ii)} $\hat{\Sigma}_m$ is complete with respect to $B_{r_1}(p_m)$, and along with $\partial B_{r_1}(p_m)$ bounds a convex set, $\hat{K}_m$;
\medskip
\myitem{(iii)} the connected component of $\Gamma_m\minter B_{r_1}(p_m)$ containing $p_m$, which we denote by $\Gamma_{m,0}$, is itself contained in $\hat{\Sigma}_m$;
\medskip
\myitem{(iv)} the outward pointing normal over $\hat{\Sigma}_m$ makes an angle of no more than $\theta/2$ with $\msf{N}^+_m$ along $\Gamma_{m,0}$; and
\medskip
\myitem{(v)} the Gaussian curvature of $\hat{\Sigma}_m$ at the point $q$ is at least $\phi_m(q)+\epsilon$.
\medskip
\noindent Moreover, we may assume that $(\hat{\Sigma}_m)_\minn$ converges towards $\hat{\Sigma}_0$.
\medskip
\noindent By Lemma \procref{LemmaBoundaryCompactnessWhenNotEmbedded}, reducing $r_1$ if necessary we may assume that, for all $m$, the connected component of the intersection of $\Sigma_m$ with $B_{r_1}(p_m)$ containing $p_m$, which we denote by $\Sigma_{m,0}$, is embedded, lies on the boundary of a convex set, $K_m$ such that $K_m\subseteq\hat{K}_m$. By compactness of the family of compact sets, there exists a convex set $K_0$ to which $(K_m)_\minn$ converges in the Haussdorf sense. The angle that the normal to $K_0$ makes with $T\hat{\Sigma}_0$ at $p_0$ is strictly less than $\pi$. Thus, for all $m$, $\Sigma_m$ is a graph over some (almost) fixed hypersurface over a uniform radius about $p$: formally, reducing $r_1$ further if necessary, for all $m$, there exists a smooth embedded hypersurface $S_m\subseteq M$ and an open subset $\Omega_m\subseteq S_m$ with smooth boundary such that:
\medskip
\myitem{(i)} $p_m\in S_m$ and $S_m$ is complete with respect to $B_{r_1}(p_m)$;
\medskip
\myitem{(ii)} the shape operator of $S_m$ vanishes at $p_m$;
\medskip
\myitem{(iii)} $\Gamma_m$ is a graph over $\partial\Omega_m$; and
\medskip
\myitem{(iv)} $\Sigma_{m,0}$ and $\hat{\Sigma}_m$ are graphs of functions $f_m$ and $\hat{f}_m$ respectively over $\Omega_m$ such that $\hat{f}_m\geqslant f_m$.
\medskip
\noindent Moreover, we may suppose that $(S_m)_\minn$ converges to $S_0$ and that $(\hat{f}_m)_\minn$ converges to $\hat{f}_0$ in the $C^\infty_\oploc$ sense. Using this construction in conjunction with Proposition $5.1$ of \cite{SmiCGC} and Proposition \procref{PropFirstOrderLowerBound}, we obtain $K_1>0$ such that, for all $m$ and for all $p\in\Gamma_m$:
$$
\|A_m(p)\| \leqslant K_1.
$$
\noindent Since the diameter of $\Sigma_m$ is uniformly bounded above, by Proposition $6.1$ of \cite{SmiCGC}, we obtain $K_2>0$ such that, for all $m$, and for all $p\in\Sigma_m$:
$$
\|A_m(p)\| \leqslant K_2.
$$
\noindent Again, using the above construction along with Theorem $1$ of \cite{CaffNirSprII}, we show that there exists $\epsilon>0$ and uniform $C^{0,\alpha}$ bounds for $(A_m)_\minn$. The Schauder estimates then yield uniform $C^k$ bounds for $(A_m)_\minn$ for all $k$. The result now follows by the Arzela-Ascoli Theorem of \cite{SmiAAT}.\qed
\medskip
\noindent Let $(\hat{\Sigma}_m)_\minn,\Sigma_0\subseteq M$ be locally strictly convex, immersed hypersurfaces in $M$ with generic boundaries such that $(\hat{\Sigma}_m)_\minn$ converges to $\Sigma_0$. Let $(\phi_m)_\minn,\phi_0:M\rightarrow]0,\infty[$ be smooth, positive functions such that $(\phi_m)_\minn$ converges to $\phi_0$ in the $C^\infty_\oploc$ sense.
\medskip
\noindent Lemma \procref{LemmaPreCompactness} can be refined to the following result:
\proclaim{Lemma \nextprocno}
\noindent Let $(\Sigma_m)_\minn$ be strictly convex smooth immersed hypersurfaces in $M$ such that, for all $m$:
\medskip
\myitem{(i)} $\Sigma_m$ is bounded by $\hat{\Sigma}_m$; and
\medskip
\myitem{(ii)} for all $p\in\Sigma_m$, the Gaussian curvature of $\Sigma_m$ at $p$ is equal to $\phi_m(p)$.
\medskip
\noindent There exists a strictly convex smooth immersed hypersurface, $\Sigma_0$ in $M$ to which $(\Sigma_m)_\minn$ subconverges. Moreover:
\medskip
\myitem{(i)} $\Sigma_0$ is bounded by $\hat{\Sigma}_0$; and
\medskip
\myitem{(ii)} for all $p\in\Sigma_0$, the Gaussian curvature of $\Sigma_0$ at $p$ is equal to $\phi_0(p)$.
\endproclaim
\proclabel{LemmaRefinedPrecompactness}
\proof Since $(\hat{\Sigma}_m)_\minn$ converges to $\hat{\Sigma}_0$, there exists $D>0$ such that, for all $m$, the diameter of $\hat{\Sigma}_m$ is bounded above by $D$. Likewise, for all $m$, $\Gamma_m:=\partial\hat{\Sigma}_m$ is locally strictly convex and, if $\msf{N}^-_m$ and $\msf{N}^+_m$ denote the convexity orientation and coorientation respectively of $\Gamma_m$, then there exists $\theta>0$ such that the angle that the outward pointing unit normal to $\hat{\Sigma}_m$ makes with $\msf{N}^+_m$ along $\Gamma_m$ is everywhere bounded below by $\theta$.
\medskip
\noindent For all $m$, let $\pi_m:\hat{\Sigma}_m\rightarrow\Sigma_m$ be the canonical projection. Since $M$ has non-positive curvature, for all $m$, $\pi_m$ is distance decreasing, and the diameter of $\Sigma_m$ is thus bounded above by $D$. Moreover, for all $m$, since $\hat{\Sigma}_m$ bounds $\Sigma_m$, the angle that the outward pointing unit normal to $\Sigma_m$ makes with $\msf{N}^+_m$ along $\Gamma_m$ is everywhere bounded below by $\theta$. It follows by Lemma \procref{LemmaPreCompactness} that there exists a strictly convex immersed hypersurface, $\Sigma_0$ towards which $(\Sigma_m)_\minn$ subconverges such that, for all $p\in\Sigma_0$, the Gaussian curvature of $\Sigma_0$ at $p$ is equal to $\phi_0(p)$. By Lemma \procref{LemmaClosednessOfGraphProperty}, $\hat{\Sigma}_0$ bounds $\Sigma_0$ and this completes the proof.\qed
\goodbreak
\newhead{Local Deformation}
\noindent Let $M^{n+1}$ be a Hadamard manifold. Let $(\hat{\Sigma}_t)_{t\in[0,1]}$ be a smooth family of locally convex immersed hypersurfaces in $M$ with generic boundary. For all $t$, denote $\Gamma_t=\partial\hat{\Sigma}_t$. Let $\epsilon>0$ and let $(\phi_t)_{t\in[0,1]}\in C^\infty(M,]0,\infty[)$ be a smooth family such that, for all $t$, the Gaussian curvature of $\hat{\Sigma}_t$ is everywhere greater than $\phi_t+\epsilon$.
\headlabel{HeadLocalDeformation}
\medskip
\noindent For all $t\in[0,1]$ let $\Cal{M}_t$ be as in Section \headref{HeadImmersedSubmanifolds} and let $\Cal{N}_t$ be the family of (equivalence classes) of convex immersed hypersurfaces, $[\Sigma]$ in $M$ such that $\partial\Sigma=\partial\hat{\Sigma}_t$ and $\Sigma$ is strictly bounded by $\hat{\Sigma}_t$. By Lemma \procref{LemmaOpennessOfGraphProperty}, $\Cal{N}_t$ is an open subset of $\Cal{M}_t$ and is therefore interpreted as a smooth Banach manifold. Let $\Cal{M}$ be as in Section \headref{HeadImmersedSubmanifolds} and let $\Cal{N}$ be the family of all pairs $(t,[\Sigma])$ where $t\in[0,1]$ and $[\Sigma]\in\Cal{N}_t$. $\Cal{N}$ is likewise an open subset of $\Cal{M}$.
\medskip
\noindent Let $X_0\subseteq\Cal{N}$ be the set of all pairs $(t,[\Sigma])$ in $\Cal{N}$ such that the Gaussian curvature of $\Sigma$ is equal to $\phi_t$. By Lemma \procref{LemmaRefinedPrecompactness} and the Geometric Maximum Principal, $X_0$ is compact. Let $P=(t_0,[\Sigma])$ be a point in $X_0$, where $\Sigma=(i,(S,\partial S))$.
Let $(i_t)_{t\in]t_0-\epsilon,t_0+\epsilon[}$ be a smooth family of immersions such that $i_0=i$ and, for all $t$, $\Gamma_t=(i_t,\partial S)$. We define the family $(\Sigma_t)_{t\in]t_0-\epsilon,t_0+\epsilon[}$ by:
$$
\Sigma_s = (i_s,(S,\partial S)).
$$
\noindent Let $(U_P,V_P,\Phi_P)$ be the resulting graph neighbourhood of $\Cal{N}$ about $\Sigma$.
\medskip
\noindent Consider the Gauss curvature mapping $K$. This is a smooth section of $\Cal{E}$. If we identify $T_P\Cal{N}_t$ with $C_0^\infty(S)$, then its covariant derivative, $\nabla K$, defines a mapping from $C_0^\infty(S)$ to $C^\infty(S)$. By Corollary \procref{CorInvertibility}, $\nabla K$ is a second order elliptic linear differential operator. It is therefore Fredholm. Since it maps from $C^\infty_0(S)$ to $C^\infty(S)$, it is of index $0$. There therefore exists a finite dimensional vector subspace $E\subseteq C^\infty(S)$ such that if $M$ is defined by:
$$
M:E\oplus C_0^\infty(S)\rightarrow C^\infty(S);(f,\phi)\mapsto \nabla K\cdot\phi + f,
$$
\noindent then $M$ is surjective. Since $M$ differs from $\nabla K$ by a compact (in fact, finite rank) operator, it is Fredholm of index $m$, where $m$ is the dimension of $E$. Let $f_1,...,f_n$ be a basis of $E$. For $Q:=(t_Q,\Sigma_Q)\in U_P$, where $\Sigma_Q = (i_Q,(S_Q,\partial S_Q))$, let $\pi_Q:(S_Q,\partial S_Q)\rightarrow(S,\partial S)$ be the canonical projection (recall that $\Sigma_Q$ is a graph over $\Sigma_{t_Q}$). For all $i$, we define $f_{i,Q}\in C^\infty(S_Q)$ by:
$$
f_{i,Q} = f_i\circ\pi_Q.
$$
\noindent For all $i$, $Q\mapsto f_{i,Q}$ defines a section of $\Cal{E}|_{U_P}$, which we denote by $F_i$. We now define $\hat{K}_P:\Bbb{R}^m\times U_P\rightarrow\Cal{E}|_{U_P}$ by:
$$
\hat{K}_P(\sum_{i=1}^n\lambda_i e_i,(t,[\Sigma])) = K(\Sigma) + \sum_{i=1}^n\lambda_i F_i(t,[\Sigma]).
$$
\noindent By reducing $U_P$ if necessary, we may assume that $\nabla\hat{K}_P$ is Fredholm and surjective at every point of $\Bbb{R}^m\times U_P$. Since $\hat{K}_P$ is now a function over an open subset of $\Cal{M}$ (as opposed to $\Cal{M}_t$), it's derivative has index $(m+1)$.
\medskip
\noindent More generally, let $\psi:U_P\rightarrow[0,\infty[$ be a smooth function such that:
\medskip
\myitem{(i)} $\psi=1$ near $(t_0,[\Sigma])$; and
\medskip
\myitem{(ii)} the support of $\psi$ is contained in $U_P$.
\medskip
\noindent Let $U'_P\subseteq U_P$ be a neighbourhood of $(t_0,[\Sigma])$ such that $\psi=1$ over $U_P'$. We define $\Psi_P:\Bbb{R}^m\rightarrow\Gamma(\Cal{E})$ by:
$$
\Psi_P(\sum_{i=1}^n\lambda_i e_i) = \sum_{i=1}^n\lambda_i\psi F_i.
$$
\noindent By compactness of $X_0$, there exist finitely many points $P_1,...,P_n\in X_0$ such that:
$$
X_0\subseteq\munion_{i=1}^n U_{P_i}'=:\Omega.
$$
\noindent Denote $m=m_1+...+m_n$ and define $\Psi:\Bbb{R}^m\rightarrow\Gamma(\Cal{E})$ by:
$$
\Psi = \Psi_{P_1}\oplus...\oplus\Psi_{P_n}.
$$
\noindent Define $\hat{K}:\Bbb{R}^m\times\Cal{N}\rightarrow\Cal{E}$ by:
$$
\hat{K}(v,(t,[\Sigma])) = K([\Sigma]) + \Psi(v).
$$
\noindent For $v\in\Bbb{R}^m$, define $X_v$ by:
$$
X_v = \left\{(t,[\Sigma])\in\Cal{N}\text{ s.t. }\hat{K}(v,(t,[\Sigma])) = \phi_t\right\}.
$$
\proclaim{Proposition \nextprocno}
\noindent There exists $r>0$ such that:
\medskip
\myitem{(i)} for $\|v\|<r$, $X_v$ is compact; and
\medskip
\myitem{(ii)} for $\|v\|<r$, $X_v\subseteq\Omega$.
\endproclaim
\proclabel{PropTheSetsAreCompactAndConverge}
\proof $(i)$ Let $(t_m,[\Sigma_m])_\minn$ be a sequence in $X_v$. Let $(\Sigma'_m)_\minn$ be a sequence of smooth, immersed, compact hypersurfaces in $M$ such that, for all $m$, $\Sigma_m$ is a graph over $\Sigma'_m$. Suppose, moreover, that $(\Sigma'_m)_\minn$ converges to $\Sigma_0'$. For all $m\in\Bbb{N}\munion\left\{0\right\}$, choose $f_m\in C^\infty(\Sigma'_m)$ and suppose that $(f_m)_\minn$ converges in the $C^\infty$ sense to $f_0$. For all $m$, let $\pi_m$ be the canonical projection onto $\Sigma'_m$. With small modifications, Lemma \procref{LemmaRefinedPrecompactness} adapts to the case where $\phi_m = f_m\circ\pi_m$ for all $m$, and likewise to the case where $\phi_m$ is a finite linear combination of such functions. It follows that the closure of $X_v$ in $\Cal{M}$ is relatively compact.
\medskip
\noindent Let $(t,[\Sigma])$ be a limit point of $X_v$. By Lemma \procref{LemmaOpennessOfGraphProperty}, $\Sigma$ is bounded by $\hat{\Sigma}_t$. Suppose that $\Sigma\notin\Cal{N}_t$. Then $\hat{\Sigma}_t$ does not strictly bound $\Sigma$, and $\Sigma$ is thus an interior tangent to $\hat{\Sigma}_t$ at some point, $p$, say (possibly in $\partial\hat{\Sigma}$). However, for $v$ sufficiently small, $\|\Psi(v)\|\leqslant \epsilon$ and so the Gaussian curvature of $\hat{\Sigma}_t$ at $p$ is strictly greater than that of $\Sigma$ at $p$. This contradicts the Geometric Maximum Principal (see, for example, \cite{SmiCGC}). There thus exists $r>0$ such that for $\|v\|<r$, the closure of $X_v$ is contained in $\Cal{N}$ and so $X_v$ is compact. $(i)$ follows.
\medskip
\noindent $(ii)$ Suppose the contrary. There exists $(v_n)_\ninn$ which converges to $0$ and $(t_n,[\Sigma_n])_\ninn$ such that, for all $n$:
$$
(t_n,[\Sigma_n])\in X_{v_n},\qquad (t_n,[\Sigma_n])\notin\Omega.
$$
\noindent As in the previous paragraph, by Lemma \procref{LemmaRefinedPrecompactness}, $(t_n,[\Sigma_n])_\ninn$ subconverges to $(t_0,[\Sigma_0])\in X_0$. Thus, for sufficiently large $n$, $(t_n,(\Sigma_n))_\ninn\in\Omega$, which is absurd. $(ii)$ follows, and this completes the proof.\qed
\medskip
\noindent Define $X\subseteq\Bbb{R}^m\times\Omega$ by:
$$
X = \left\{(v,(t,[\Sigma]))\in\Bbb{R}^m\times\Omega\text{ s.t. }\hat{K}(v,(t,[\Sigma])) = \phi_t\right\}.
$$
\proclaim{Proposition \nextprocno}
\noindent $X$ is an $(m+1)$-dimensional smooth, embedded submanifold of $\Bbb{R}^m\times\Omega$.
\endproclaim
\proof By construction, $\hat{K}$ is everywhere Fredholm of index $(m+1)$ and surjective. The result now follows by the Implicit Function Theorem for Banach manifolds.\qed
\proclaim{Proposition \nextprocno}
\noindent There exists $(v_n)_\ninn\in\Bbb{R}^m$ such that:
\medskip
\myitem{(i)} $(v_n)_\ninn$ converges to $0$;
\medskip
\myitem{(ii)} for all $n$, $X_{v_n}$ is a (potentially empty) $1$-dimensional, smooth, compact, embedded submanifold of $\Omega$; and
\medskip
\myitem{(iii)} $\partial X_{v_n}\subseteq \Cal{N}_0\munion\Cal{N}_1$.
\endproclaim
\proclabel{PropSmoothApproximations}
\proof Let $\pi:\Bbb{R}^m\times\Omega\rightarrow\Bbb{R}^m$ be projection onto the first factor. Let $\pi_X$ be the restriction of $\pi$ to $X$. By Sard's Lemma, the set of critical values of $\pi_X$ has Lebesgue measure $0$. Let $(v_n)_\ninn\in\Bbb{R}^m$ be a sequence of non-critical values of $\pi_X$ converging to $0$. By the Submersion Theorem, for all $n$, $X_{v_n}$ is a $1$-dimensional, smooth, embedded submanifold of $X$ and therefore of $\Omega$. By Proposition \procref{PropTheSetsAreCompactAndConverge} we may suppose moreover that, for all $n$, $X_{v_n}$ is compact. $(i)$ and $(ii)$ follow. For all $n$, the end points of $X_{v_n}$ lie in the (manifold) boundary of $X$. Since this is contained in $\Cal{N}_0\munion\Cal{N}_1$, $(iii)$ follows. This completes the proof.\qed
\goodbreak
\newhead{Local and Global Rigidity}
\noindent Let $M^{n+1}$ be an $(n+1)$-dimensional Hadamard manifold. Let $\hat{\Sigma}\subseteq M$ be a convex immersed hypersurface. Choose $\phi\in C^\infty(M)$. Let $\Sigma=(i,(S,\partial S))$ be another convex immersed hypersurface. We say that $\Sigma$ is a {\bf solution} to the problem $(\hat{\Sigma},\phi)$ if and only if:
\medskip
\myitem{(i)} $\partial\Sigma=\partial\hat{\Sigma}$;
\medskip
\myitem{(ii)} $\Sigma$ is bounded by $\hat{\Sigma}$; and
\medskip
\myitem{(iii)} for all $p\in S$, the Gaussian curvature of $\Sigma$ at $p$ is equal to $(\phi\circ i)(p)$.
\headlabel{HeadLocalAndGlobalRigidity}
\proclaim{Definition \nextprocno}
\myitem{(i)} We say that $(\hat{\Sigma},\phi)$ is locally rigid if and only if, for all solutions, $\Sigma$ to $(\hat{\Sigma},\phi)$, the linearisation, $DK$, of the Gauss Curvature Operator, $K$, over $\Sigma$ is invertible.
\medskip
\myitem{(ii)} We say that $(\hat{\Sigma},\phi)$ is globaly rigid if and only if there exists at most one solution, $\Sigma$ to $(\hat{\Sigma},\phi)$.
\endproclaim
\proclabel{DefnLocalAndGlobalRigidity}
\noindent We recall the following properties of local and global rigidity:
\proclaim{Proposition \nextprocno}
\myitem{(i)} If $(\hat{\Sigma},\phi)$ is locally rigid, then $(\hat{\Sigma},\phi')$ is also locally rigid for all $\phi'$ sufficiently close to $\phi$.
\medskip
\myitem{(ii)} If $(\hat{\Sigma},\phi)$ is locally and globally rigid, then $(\hat{\Sigma},\phi')$ is globally rigid for all $\phi'$ sufficiently close to $\phi$.
\endproclaim
\proclabel{PropPropertiesOfRigidity}
\proof See \cite{SmiCGC}.\qed
\medskip
\noindent Now let $(\hat{\Sigma}_t)_{t\in[0,1]}$ be a smooth family of locally strictly convex, immersed hypersurfaces in $M$ with generic boundaries. Let $\epsilon>0$ and let $(\phi_t)_{t\in[0,1]}\in C^\infty(M,]0,\infty[)$ be a smooth family of smooth, positive functions such that, for all $t$, the Gaussian curvature of $\Sigma_t$ at any point $p$ is no less than $\phi_t(p)+\epsilon$. Using local and global rigidity, we obtain existence:
\proclaim{Lemma \nextprocno}
\noindent Suppose that $(\hat{\Sigma}_0,\phi_0)$ is both locally and globally rigid. If there exists a solution $\Sigma_0$ to $(\hat{\Sigma}_0,\phi_0)$, then there exists a solution to $(\hat{\Sigma}_1,\phi_1)$.
\endproclaim
\proclabel{LemmaExistence}
\remark It follows that proving existence of solutions for a given problem reduces to showing the existence of a smooth isotopy by locally strictly convex immersions to a locally and globally rigid problem for which solutions are known to exist.
\medskip
\proof Let $\Cal{N}$, $m\in\Bbb{N}$ and $\Psi:\Bbb{R}^m\rightarrow\Gamma(\Cal{E})$ be as in Section \headref{HeadLocalDeformation} and, for all $v\in\Bbb{R}^m$, define $X_v\subseteq\Cal{N}$ by:
$$
X_v = \left\{(t,[\Sigma])\in\Cal{N}\text{ s.t. }K([\Sigma]) + \Psi(v) = \psi_t\right\}.
$$
\noindent Let $(v_n)_\ninn\subseteq\Bbb{R}^m$ be as in Proposition \procref{PropSmoothApproximations}. Since $(\hat{\Sigma}_0,\phi_0)$ is locally rigid, there exists $N>0$ such that, for all $n\geqslant N$, $X_{v_n}\minter\Cal{N}_0$ is non-empty, and thus, in particular, $X_{v_n}$ is non-empty. Since $(\hat{\Sigma}_0,\phi_0)$ is also globally rigid, it follows by Proposition \procref{PropPropertiesOfRigidity} that, for sufficiently large $n$, $\Psi(v_n) + \phi_0$ is too, and therefore that $X_{v_n}\minter\Cal{N}_0$ consists of a single point.
\medskip
\noindent Let $\pi:\Cal{N}\rightarrow[0,1]$ be the canonical projection. For all $n\geqslant N$, $X_{v_n}$ is a smooth, embedded, compact, $1$-dimensional submanifold of $\Cal{N}$. It is thus homeomorphic, either to a compact interval or to a circle. By local and global rigidity, the restriction of $\pi$ to $X_{v_n}$ is a local diffeomorphism near the unique point lying in $\pi^{-1}(\left\{0\right\})$. It follows that $X_{v_n}$ has non-trivial (manifold) boundary, and is therefore not a circle. It is thus a compact interval. By Proposition \procref{PropSmoothApproximations}, the endpoints of $X_{v_n}$ lie in $\Cal{N}_0\munion\Cal{N}_1$. By global rigidity, only one endpoint of $X_{v_n}$ lies in $\Cal{N}_0$, and the other therefore lies in $\Cal{N}_1$.
\medskip
\noindent For all $n$, let $\Sigma_n$ be such that $(1,[\Sigma_n])$ is the unique endpoint of $X_{v_n}$ in $\Cal{N}_1$. By Lemma \procref{LemmaRefinedPrecompactness}, there exists $\Sigma_0$ to which $(\Sigma_n)_\ninn$ subconverges and $\Sigma_0$ is a solution of $(\hat{\Sigma}_1,\psi_1)$. This completes the proof.\qed
\medskip
\noindent Lemma \procref{LemmaExistence} may be easily adapted to treat the case where the metric of the underlying manifold also varies, and we obtain Theorem \procref{ThmExistenceII}:
\medskip
{\bf\noindent Proof of Theorem \procref{ThmExistenceII}:\ }Let $(\hat{\Sigma}_t)_{t\in[0,1]}$ be an isotopy by convex, immersed hypersurfaces such that $\hat{\Sigma}_0=\hat{\Sigma}$ and $\hat{\Sigma}_1$ is a finite covering of $\Omega$. For ease of presentation, we will assume that the covering is of order one: the general case is almost identical.
Let $p\in K$ be an interior point. Let $d_0,d_1:M\rightarrow\Bbb{R}$ be given by:
$$
d_0(x) = d(x,K),\qquad d_1(x) = d(x,p).
$$
\noindent Both $d_0$ and $d_1$ are smooth outside $K$. For $t\in[0,1]$, define $d_t$ by:
$$
d_t = td_1 + (1-t)d_0.
$$
\noindent Trivially, $\partial K$ is isotopic by smooth convex immersions to $d_0^{-1}(\left\{r\right\})$ for all $r\geqslant 0$. Choose $r_0$ such that $K\subseteq B_{r_0}(p)$. For all $t$, $d_t^{-1}(\left\{r_0\right\})$ is a convex, embedded hypersurface and we thus obtain an isotopy by smooth convex immersions between $d_0^{-1}(\left\{r_0\right\})$ and $d_1^{-1}(\left\{r_0\right\})$. We may thus define $(\hat{\Sigma}_t)_{t\in[1,2]}$ such that $\hat{\Sigma}_2$ is a geodesic sphere with a finite number of open sets removed. Let $g$ be the Riemannian metric on $M$. Define $(g_t)_{t\in[0,2]}$ such that $g_t=g$ for all $t$.
\medskip
\noindent We may assume that $\hat{\Sigma}_2$ is as small as we wish. Define $(\hat{\Sigma}_t)_{t\in[2,3]}$ and $(g_t)_{t\in[2,3]}$ such that:
\medskip
\myitem{(i)} $g_2=g$;
\medskip
\myitem{(ii)} $g_3$ is complete with constant curvature equal to $1$;
\medskip
\myitem{(iii)} for all $t$, $\hat{\Sigma}_t$ is a geodesic sphere with respect to $g_t$ with a finite number of open sets removed.
\medskip
\noindent Define $(\hat{\Sigma}_t)_{t\in[3,4]}$ and $(g_t)_{t\in[3,4]}$ such that:
\medskip
\myitem{(i)} for all $t$, $g_t=g_3$ is the complete hyperbolic metric;
\medskip
\myitem{(ii)} for all $t$, $(\hat{\Sigma}_t)$ is a geodesic sphere with a finite number of open sets removed; and
\medskip
\myitem{(iii)} $\hat{\Sigma}_4$ is a horosphere with a finite number of open sets (including a neighbourhood of the infinite point) removed.
\medskip
\noindent Let $(\psi_t)_{t\in[0,4]}\in C^\infty(M)$ be a smooth family of smooth, positive valued functions such that:
\medskip
\myitem{(i)} $\psi_0=\psi$;
\medskip
\myitem{(ii)} for all $t$ and for all $p\in\hat{\Sigma}_t$, the Gaussian curvature of $\hat{\Sigma}_t$ at $p$ is greater than $\psi_t(p)$; and
\medskip
\myitem{(iii)} $\psi_4$ is constant and equal to $1-\delta$ for some $\delta<1$.
\medskip
\noindent The problem $(\Sigma_4,\psi_4)$ in $(M,g_4)=\Bbb{H}^{n+1}$ is locally and globally rigid and has a non-trivial solution (see \cite{SmiCGC}). By Proposition \procref{PropDensityOfGenericity}, this isotopy by locally strictly convex, immersed hypersurfaces may be deformed to an isotopy by locally strictly convex, immersed hypersurfaces whose boundaries are generic. Existence therefore follows by (an appropriately modified version of) Lemma \procref{LemmaExistence}, and this completes the proof.\qed
\inappendicestrue
\global\headno=0
\goodbreak
\newhead{Bibliography}
{\leftskip = 5ex \parindent = -5ex
\leavevmode\hbox to 4ex{\hfil \cite{HarveyLawsonI}}\hskip 1ex{Harvey F. R., Lawson H. B. Jr., Dirichlet Duality and the Nonlinear Dirichlet Problem, {\sl Comm. Pure Appl. Math.} {\bf 62} (2009), no. 3, 396--443}
\medskip
\leavevmode\hbox to 4ex{\hfil \cite{CaffNirSprI}}\hskip 1ex{Caffarelli L., Nirenberg L., Spruck J., The Dirichlet problem for nonlinear second-Order elliptic equations. I. Monge Amp\`ere equation, {\sl Comm. Pure Appl. Math.} {\bf 37} (1984), no. 3, 369--402}
\medskip
\leavevmode\hbox to 4ex{\hfil \cite{CaffNirSprII}}\hskip 1ex{Caffarelli L., Kohn J. J., Nirenberg L., Spruck J., The Dirichlet problem for nonlinear second-order elliptic equations. II. Complex Monge Amp\`ere, and uniformly elliptic, equations, {\sl Comm. Pure Appl. Math.} {\bf 38} (1985), no. 2, 209--252}
\medskip
\leavevmode\hbox to 4ex{\hfil \cite{CaffNirSprV}}\hskip 1ex{Caffarelli L., Nirenberg L., Spruck J., Nonlinear second-order elliptic equations. V. The Dirichlet problem for Weingarten hypersurfaces, {\sl Comm. Pure Appl. Math.} {\bf 41} (1988), no. 1, 47--70}
\medskip
\leavevmode\hbox to 4ex{\hfil \cite{GuanSpruckO}}\hskip 1ex{Guan B., Spruck J., Boundary value problems on $\Bbb{S}^n$ for surfaces of constant Gauss curvature, {\sl Ann. of Math.} {\bf 138} (1993), 601--624}
\medskip
\leavevmode\hbox to 4ex{\hfil \cite{GuanSpruckI}}\hskip 1ex{Guan B., Spruck J., The existence of hypersurfaces of constant Gauss curvature with prescribed boundary, {\sl J. Differential Geom.} {\bf 62} (2002), no. 2, 259--287}
\medskip
\leavevmode\hbox to 4ex{\hfil \cite{GuillPoll}}\hskip 1ex{Guillemin V., Pollack A., {\sl Differential Topology\/}, Prentice-Hall, Englewood Cliffs, N.J., (1974)}
\medskip
\leavevmode\hbox to 4ex{\hfil \cite{Klingenberg}}\hskip 1ex{Klingenberg W., {\sl Lectures on closed geodesics\/}, Grundlehren der Mathematischen Wissenschaften, {\bf 230}, Springer-Verlag, Berlin-New York, (1978)}
\medskip
\leavevmode\hbox to 4ex{\hfil \cite{LabI}}\hskip 1ex{Labourie F., Un lemme de Morse pour les surfaces convexes (French), {\sl Invent. Math.} {\bf 141} (2000), no. 2, 239--297}
\medskip
\leavevmode\hbox to 4ex{\hfil \cite{RosSpruck}}\hskip 1ex{Rosenberg H., Spruck J., On the existence of locally convex hypersurfaces of constant Gauss curvature in hyperbolic space, {\sl J. Differential Geom.} {\bf 40} (1994), no. 2, 379--409}
\medskip
\leavevmode\hbox to 4ex{\hfil \cite{SmiCGC}}\hskip 1ex{Smith G., Compactness for immersions of prescribed Gaussian curvature I - analytic aspects, {\sl in preparation}}
\medskip
\leavevmode\hbox to 4ex{\hfil \cite{SmiNLD}}\hskip 1ex{Smith G., The Non-Linear Dirichlet Problem in Hadamard Manifolds,\break arXiv:0908.3590}
\medskip
\leavevmode\hbox to 4ex{\hfil \cite{SmiAAT}}\hskip 1ex{Smith G., An Arzela-Ascoli Theorem for immersed submanifolds, {\sl Ann. Fac. Sci. Toulouse Math.}, {\bf 16}, no. 4, (2007), 817--866}
\medskip
\leavevmode\hbox to 4ex{\hfil \cite{SmiGPP}}\hskip 1ex{Smith G., The non-linear Plateau problem in non-positively curved manifolds,\break arXiv:1004.0374}
\medskip
\leavevmode\hbox to 4ex{\hfil \cite{Spruck}}\hskip 1ex{Spruck J., Fully nonlinear elliptic equations and applications to geometry, {\sl Proceedings of the International Congress of Mathematicians}, (Z\"urich, 1994), 1145--1152, Birkh\"a user, Basel, 1995.}
\medskip
\leavevmode\hbox to 4ex{\hfil \cite{TrudWang}}\hskip 1ex{Trudinger N. S., Wang X., On locally locally convex hypersurfaces with boundary, {\sl J. Reine Angew. Math.} {\bf 551} (2002), 11--32}
\medskip
\leavevmode\hbox to 4ex{\hfil \cite{White}}\hskip 1ex{White B., The space of minimal submanifolds for varying Riemannian metrics, {\sl Indiana Univ. Math. J.} {\bf 40} (1991),  no. 1, 161--200}
\par}
\enddocument